\documentclass[11pt,a4paper]{amsart}[2000]
\usepackage{geometry}
\usepackage{tikz}
\usepackage{amsmath}
\usepackage{amssymb}
\usepackage{amsthm}
\usepackage[hang]{footmisc}
\usepackage[all]{xypic}
\usepackage[colorlinks=true,linkcolor=blue,citecolor=blue]{hyperref}
\usepackage{cleveref}
\usepackage{color}
\usepackage{mathtools}
\usepackage{enumitem}
\usepackage{stackengine}
\usepackage{caption}
\usetikzlibrary{arrows.meta, decorations.markings}

\title[Almost elementary models]{Almost elementary groupoid models for $C^*$-algebras}

\thanks{}

\allowdisplaybreaks

\theoremstyle{plain}

\newtheorem{Thm}{Theorem}[section]

\theoremstyle{definition}

\newtheorem{Exl}[Thm]{Example}

\theoremstyle{plain}
\newtheorem{thm}[Thm]{Theorem}
\newtheorem{lem}[Thm]{Lemma}
\newtheorem{cor}[Thm]{Corollary}
\newtheorem{prop}[Thm]{Proposition}

\newtheorem{thmx}{Theorem}
\newtheorem{corx}[thmx]{Corollary}
\newtheorem{conjx}{Conjecture}
\newtheorem{quesx}[conjx]{Question}

\theoremstyle{definition}
\newtheorem{defn}[Thm]{Definition}

\newtheorem{rmk}[Thm]{Remark}

\newcommand{\A}[0]{\mathbb{A}}
\newcommand{\B}{B}
\newcommand{\J}{J}
\newcommand{\K}{\mathcal{K}}
\newcommand{\D}{D}
\newcommand{\Ch}{D}
\newcommand{\Zh}{\mathcal{Z}}
\newcommand{\E}{E}
\newcommand{\Oh}{\mathcal{O}}

\newcommand{\T}{{\mathbb T}}

\newcommand{\N}{{\mathbb N}}
\newcommand{\Z}{{\mathbb Z}}
\newcommand{\C}{{\mathbb C}}
\newcommand{\Q}{{\mathbb Q}}

\newcommand{\aut}{\mathrm{Aut}}
\newcommand{\supp}{\mathrm{supp}}

\newcommand{\eps}{\varepsilon}
\numberwithin{equation}{section}

\newcommand{\id}{\mathrm{id}}

\newcommand{\halpha}{\widehat{\alpha}}
\newcommand{\calpha}{\widehat{\alpha}}

\newcommand{\tih}{\widetilde {h}}

\newcommand\set[1]{\left\{#1\right\}}  
\newcommand\mset[1]{\left\{\!\!\left\{#1\right\}\!\!\right\}}

\newcommand{\CA}[0]{\mathcal{A}} \newcommand{\CB}[0]{\mathcal{B}}
\newcommand{\CC}[0]{\mathcal{C}} \newcommand{\CD}[0]{\mathcal{D}}
 
\newcommand{\CG}[0]{\mathcal{G}} \newcommand{\CH}[0]{\mathcal{H}}
 
\newcommand{\CK}[0]{\mathcal{K}} \newcommand{\CL}[0]{\mathcal{L}}
 
\newcommand{\CO}[0]{\mathcal{O}} 
\newcommand{\CQ}[0]{\mathcal{Q}} \newcommand{\CR}[0]{\mathcal{R}}
 \newcommand{\CT}[0]{\mathcal{T}}
\newcommand{\CU}[0]{\mathcal{U}} \newcommand{\CV}[0]{\mathcal{V}}
 
 \newcommand{\CZ}[0]{\mathcal{Z}}

\newcommand{\Ra}[0]{\Rightarrow}
\newcommand{\La}[0]{\Leftarrow}
\newcommand{\LRa}[0]{\Leftrightarrow}

\newcommand{\quer}[0]{\overline}
\newcommand{\eins}[0]{\mathbf{1}}			
\newcommand{\diag}[0]{\operatorname{diag}}
\newcommand{\ad}[0]{\operatorname{Ad}}
\newcommand{\ev}[0]{\operatorname{ev}}
\newcommand{\fin}[0]{{\subseteq\!\!\!\subseteq}}
\newcommand{\diam}[0]{\operatorname{diam}}
\newcommand{\Hom}[0]{\operatorname{Hom}}
\newcommand{\dst}[0]{\displaystyle}
\newcommand{\spp}[0]{\operatorname{supp}}
\newcommand{\lsc}[0]{\operatorname{Lsc}}
\newcommand{\del}[0]{\partial}
\newcommand{\GU}[0]{\CG^{(0)}}

\newcommand{\Ell}[0]{\operatorname{Ell}}

\theoremstyle{definition}

\numberwithin{equation}{Thm}

\begin{document}
	\global\long\def\floorstar#1{\lfloor#1\rfloor}
	\global\long\def\ceilstar#1{\lceil#1\rceil}	
	
	\global\long\def\B{B}
	\global\long\def\A{A}
	\global\long\def\J{J}
	\global\long\def\K{\mathcal{K}}
	\global\long\def\D{D}
	\global\long\def\Ch{D}
	\global\long\def\Zh{\mathcal{Z}}
	\global\long\def\E{E}
	\global\long\def\Oh{\mathcal{O}}

	\global\long\def\T{{\mathbb{T}}}
	\global\long\def\BR{{\mathbb{R}}}
	\global\long\def\N{{\mathbb{N}}}
	\global\long\def\Z{{\mathbb{Z}}}
	\global\long\def\C{{\mathbb{C}}}
	\global\long\def\Q{{\mathbb{Q}}}

	\global\long\def\aut{\mathrm{Aut}}
	\global\long\def\supp{\mathrm{supp}}

	\global\long\def\eps{\varepsilon}

	\global\long\def\id{\mathrm{id}}

	\global\long\def\halpha{\widehat{\alpha}}
	\global\long\def\calpha{\widehat{\alpha}}

	\global\long\def\tih{\widetilde{h}}

	\global\long\def\opFol{\operatorname{F{\o}l}}

	\global\long\def\opRange{\operatorname{Range}}

	\global\long\def\opIso{\operatorname{Iso}}

	\global\long\def\dimnuc{\dim_{\operatorname{nuc}}}

	\global\long\def\set#1{\left\{  #1\right\}  }

	
	\global\long\def\mset#1{\left\{  \!\!\left\{  #1\right\}  \!\!\right\}  }

	\global\long\def\Ra{\Rightarrow}
	\global\long\def\La{\Leftarrow}
	\global\long\def\LRa{\Leftrightarrow}

	\global\long\def\quer{\overline{}}
	\global\long\def\eins{\mathbf{1}}
	\global\long\def\diag{\operatorname{diag}}
	\global\long\def\ad{\operatorname{Ad}}
	\global\long\def\ev{\operatorname{ev}}
	\global\long\def\fin{{\subseteq\!\!\!\subseteq}}
	\global\long\def\diam{\operatorname{diam}}
	\global\long\def\Hom{\operatorname{Hom}}
	\global\long\def\dst{{\displaystyle }}
	\global\long\def\spp{\operatorname{supp}}
	\global\long\def\spo{\operatorname{supp}_{o}}
	\global\long\def\del{\partial}
	\global\long\def\lsc{\operatorname{Lsc}}
	\global\long\def\GU{\CG^{(0)}}
	\global\long\def\HU{\CH^{(0)}}
	\global\long\def\AU{\CA^{(0)}}
	\global\long\def\BU{\CB^{(0)}}
	\global\long\def\CUU{\CC^{(0)}}
	\global\long\def\DU{\CD^{(0)}}
	\global\long\def\QU{\CQ^{(0)}}
	\global\long\def\TU{\CT^{(0)}}
	\global\long\def\CUUU{\CC'{}^{(0)}}
	
	\global\long\def\bA{\mathbb{A}}
	\global\long\def\AUl{(\CA^{l})^{(0)}}
	\global\long\def\BUl{(B^{l})^{(0)}}
	\global\long\def\HUp{(\CH^{p})^{(0)}}
	\global\long\def\sym{\operatorname{Sym}}
	
	\global\long\def\properlength{proper}

	\global\long\def\interior#1{#1^{\operatorname{o}}}
	
	\author{Xin Ma}

	\address{X. Ma:  Institute for Advanced Study in Mathematics, Harbin Institute of Technology, Harbin, China, 150001}
  \email{xma17@hit.edu.cn}

\author{Jianchao Wu}

\address{J. Wu: Shanghai Center for Mathematical Sciences, Fudan University,
Shanghai, 200438, China}

\email{jianchao\_wu@fudan.edu.cn}	

\date{\today}

	\begin{abstract}
		The notion of almost elementariness for a locally compact Hausdorff \'{e}tale groupoid $\CG$ with a compact unit space was introduced by the authors as a sufficient condition ensuring the reduced groupoid $C^*$-algebra $C^*_r(\CG)$ is (tracially) $\CZ$-stable and thus classifiable under additional natural assumptions. 
		In this paper, we explore the converse direction and show that many groupoids in the literature serving as models for classifiable $C^*$-algebras are almost elementary. 
		In particular, for a large class $\CC$ of Elliott invariants and a $C^*$-algebra $A$ with $\operatorname{Ell}(A)\in \CC$, we show that $A$ is classifiable if and only if $A$ possesses a minimal, effective, amenable, second countable, almost elementary groupoid model, which leads to a groupoid-theoretic characterization of classifiability of $C^*$-algebras with certain Elliott invariants. 
        In addition, we 
        demonstrate obstructions to obtaining a transformation groupoid model for the Jiang-Su algebra $\CZ$. 
  \end{abstract}
	\maketitle

 \section{Introduction}
 Recently, spurred by the rapid progress in the theory of classification and structure of $C^*$-algebras,  
 there has been a growing recognition of the deep connections between the study of $C^*$-algebras and that of topological groupoids. 
 A large part of these recent developments build on top of the fulfillment of the main goal in the Elliott classification program, that is, a classification theorem which states that there is a class of so-called \emph{classifiable $C^*$-algebras} wherein any two members are $*$-isomorphic if and only if they have the same \emph{Elliott invariant} (see, e.g., \cite{EGLN}, \cite{GLN2020a}, \cite{GLN2020b} \cite{TWW}, \cite{CETWW} and \cite{CGSTW}). Here the Elliott invariant comprises of the (ordered) K-theory of the $C^*$-algebras as well as the tracial information associated with them.  
 
 Since the discovery of non-classifiable simple separable nuclear $C^*$-algebras , 
 a significant aspect of the Elliott classification program is concerned with characterizing the class of classifiable $C^*$-algebras. These efforts result in a rich theory of regularity properties for simple separable nuclear $C^*$-algebras, with the centerpiece being the Toms-Winter conjecture. 
 Using the successful verification of the first half of this conjecture (between finite nuclear dimension and $\mathcal{Z}$-stability) \cite{CETWW}, one can characterize a classifiable $C^*$-algebra $A$ as: 
 \begin{itemize}
 	\item[]  a simple separable nuclear $C^*$-algebra that is
 	\begin{enumerate}
 		\item \emph{$\mathcal{Z}$-stable}, i.e., $A$ tensorially absorbs the Jiang-Su algebra $\mathcal{Z}$, and 
 		\item in the \emph{UCT class}, i.e., $A$ satisfies the hypothesis of the universal coefficient theorem for $KK$-theory. 
 	\end{enumerate}
 \end{itemize}
 While it has long been known that the condition of $\mathcal{Z}$-stability is necessary since taking tensor products with $\mathcal{Z}$ does not change the Elliott invariant \cite{JiangSu}, the necessity of the UCT condition remains largely a mystery and is known as (or more precisely, equivalent to) the \emph{UCT problem} (see, e.g., \cite{BBWW}). 
 
 However, combining some recent results regarding Cartan subalgebras, 
 one may bypass the UCT problem and equivalently characterize a classifiable $C^*$-algebra $A$ as:
 \begin{itemize}
 	\item[] a twisted groupoid $C^*$-algebra of an amenable minimal effective second countable locally compact Hausdorff étale groupoid so that $A$ is $\mathcal{Z}$-stable. 
 \end{itemize} 
 Indeed, on the one hand, Li \cite{Li} showed that every classifiable $C^*$-algebra has a Cartan subalgebra, which, by the work of Renault \cite{Renault}, means that it is a twisted groupoid $C^*$-algebra of an effective locally compact Hausdorff étale groupoid (see also \cite{C-F-H}, \cite{DPS1}, \cite{DPS2}, and \cite{Sp} for various other constructions of groupoid models); on the other hand, by extending the method of Tu \cite{Tu}, Barlak and Li \cite{BarlakLi} showed that twisted groupoid $C^*$-algebras of amenable locally compact Hausdorff étale groupoids are in the UCT class. Finally, by a number of more classical results linking properties of  $C^*$-algebras and those of groupoids \cite{Anantharaman-DelarocheRenault, BCFS}, we know that a twisted groupoid $C^*$-algebra is nuclear, separable, and simple if and only if the groupoid is amenable, second countable, minimal, and effective. 
 
 Once we take this groupoid perspective towards classifiability, a natural question arises: 
 \begin{quesx}
 	What condition(s) on an (amenable minimal effective second countable) locally compact Hausdorff étale groupoid guarantee its (twisted) groupoid $C^*$-algebra is $\mathcal{Z}$-stable? 
 \end{quesx}

 Research in this direction, particularly in the special case of group actions on compact metrizable spaces, has been a staple topic in the Elliott classification program long before the aforementioned recent progress. 
 While earlier positive results often place strong assumptions on the acting group, the underlying space, and the action itself, recent results have reached a much greater generality. 
 Inspired by the construction of non-classifiable crossed product $C^*$-algebras and results on a conjecture of Toms relating regularity of crossed product $C^*$-algebras and mean dimension of dynamical systems, Kerr \cite{D} adapted Matui's notion \cite{Matui} of almost finiteness for ample groupoids to the setting of actions of amenable groups on compact metric spaces and showed that almost finite free minimal actions give rise to $\mathcal{Z}$-stable (simple, separable, nuclear, and stably finite) crossed products. 
 Almost finiteness was later related to the notion of small boundary property in topological dynamics and applied to produce many positive results on the classifiability of crossed product $C^*$-algebras. 
 In the setting of infinite $C^*$-algebras, the notion of \textit{pure infiniteness} for ample groupoids, first introduced by Matui in \cite{Matui2}, was extended to general étale groupoids by the first author in \cite{Ma2022Purely}, where it was demonstrated that a purely infinite groupoid (see Definition \ref{defn: pure inf for groupoid} below) gives rise to a groupoid $C^*$-algebra that is (strongly) purely infinite and thus $\CZ$-stable under some additional natural assumptions. 

 More recently, generalizing the stably finite setting above from group actions to groupoids and unifying it with the infinite setting, the authors introduced in \cite{M-W} a new approximation property called \textit{almost elementariness} for étale groupoids with compact unit spaces (see \Cref{defn: AE} below), 
 and showed it guarantees (tracial) $\CZ$-stability of groupoid $C^*$-algebras.  
 In our framework, the divide between finite and infinite $C^*$-algebras is governed by what we call \emph{fiberwise amenability} for étale groupoids, which we introduced in \cite{MWfiberwise} along with a coarse geometric framework on the groupoids. In the special case of a group action, fiberwise amenability corresponds to the amenability of the acting group. 
 
 Hence, combining the results in \cite{M-W} and the aforementioned classical facts, we know that if $\mathcal{G}$ is an amenable minimal effective second countable locally compact Hausdorff étale groupoid that has a compact unit space and is also almost elementary, then its groupoid $C^*$-algebra is classifiable.  
 Based on the known results and examples so far, it is tempting to ask the following partial converse to this result:

\begin{quesx}\label{ques: main question}
    Is it true that 
    every classifiable $C^*$-algebra 
    possesses 
    a locally compact Hausdorff \'{e}tale, second countable, minimal, topologically amenable, almost elementary (twisted) groupoid model?
\end{quesx}

We first remark that one cannot expect all groupoid models of a classifiable $C^*$-algebra to be almost elementary. Indeed, Joseph constructed in \cite{Jos} minimal topologically free non-almost elementary dynamical systems (to be more precise, it was shown in \cite{Jos} that those dynamical systems are not almost finite; on the other hand, it was proved in \cite{M-W} that the almost finiteness is equivalent to the almost elementariness for actions by amenable groups).
It was then proved in \cite{H-W} and \cite{GGGKN}
that some of the non-almost elementary dynamical systems constructed this way still yield classifiable crossed product $C^*$-algebras. 

In this paper, we provide a positive answer to Question \ref{ques: main question} for a large class of Elliott invariants.
Denote by $\CC$ the class of $C^*$-algebras $A$ whose Elliott invariant $\Ell(A)$ coincides with that of an $C^*$-algebra in either of the following classes:
\begin{enumerate}
    \item unital AF-algebras;
    \item unital $C^*$-algebras with $\Ell(A)\cong ((G_0, 1_{G_0}), G_1),$
    where $G_0, G_1$ are countable abelian groups;
    \item unital $C^*$-algebras $A$ whose Elliott's invariant \[\Ell(A)\cong ((\Z\oplus G_0, \Z_{>0}\oplus G_0, [1]=(k, 0)), G_1, \Delta, r),\]
    in which $G_0, G_1$ are countable abelian groups, $k\in \Z_{>0}$, $\Delta$ is a finite-dimensional Choquet simplex and $r: \Delta\to S(\Z\oplus G_0)$ is defined by $\tau\mapsto ((n, g)\mapsto n/k)$;
    \item crossed product $C^*$-algebras $C(Y)\rtimes_rG$ arising from minimal free actions by amenable groups $G$ on a compact Hausdorff spaces $Y$ with $\dim(Y)<\infty$.
    \end{enumerate}

Now, we state our main theorem as a combination of \Cref{thm: groupoid model stably finite} and \Cref{thm: groupoid model for kirchberg}.

\begin{thmx}\label{thm: main}
		Let $A$ be a unital $C^*$-algebra in the class $\CC$ above. Then $A$ is classifiable if and only if $A$ has a locally compact, Hausdorff, \'{e}tale, minimal, second countable, topological amenable, almost elementary groupoid model where the unit space is compact, i.e., there exists such a groupoid $\CG$ with a compact unit space such that $A\cong C^*_r(\CG)$.
	\end{thmx}
For a minimal free action of an amenable  on a compact metrizable space of finite covering dimension, it is still an open question asking whether the classifiability of its crossed product is characterized by the almost finiteness of the action? Theorem \ref{thm: main} shows that at least the classifiability is characterized by the almost elementariness of one of its groupoid model.

Moreover, since the above class covers all strongly self-absorbing $C^*$-algebras in the UCT class, namely $\CZ$, $\CO_2$, $\CO_\infty$, $\operatorname{UHF}$ and $\operatorname{UHF}\otimes \CO_\infty$, where $\operatorname{UHF}$ denotes an arbitrary UHF-algebra of infinite type, we obtain the following corollary.

\begin{corx}
    Every strongly self-absorbing $C^*$-algebra in the UCT class 
    has a locally compact, Hausdorff, \'{e}tale, minimal, second countable, topological amenable, almost elementary groupoid model.
    \end{corx}

These results suggests that almost elementariness may, to a certain extent, be considered as a groupoid analog of $\mathcal{Z}$-stability. 
The method employed to establish \Cref{thm: main} involves investigating groupoids in \cite{DPS1}, \cite{DPS2}, \cite{Sp}, and \cite{C-F-H}, and proving that these groupoids satisfy the property of being almost elementary as defined in Definition \ref{defn: AE}. Notably, the groupoids in \cite{DPS1} and \cite{DPS2} can be expressed as certain $\Z$-partial dynamical systems. By utilizing the concept of almost finiteness for global actions, as discussed in \cite{D}, the almost elementariness of the partial transformation groupoids associated with these partial actions is established.

In the case of purely infinite $C^*$-algebras, even stronger conclusions can be drawn, extending beyond the unital cases. Specifically, it is demonstrated that all groupoids constructed from a combination of $k$-graphs in \cite{Sp} and \cite{C-F-H} are purely infinite in  Definition \ref{defn: pure inf for groupoid}. Furthermore, the equivalence between pure infiniteness and almost elementariness is established in cases where the unit space is compact and there are no invariant measures, as shown in Proposition \ref{prop: pure inf equal AE}.
It follows that all Kirchberg algebras satisfying the UCT have a purely infinite groupoid model, as stated in \Cref{thm: groupoid model for kirchberg}.

We conclude the introduction with an exploratory discussion of the following related question. This one is also recorded as Problem XLV in \cite{STW2025}.
\begin{quesx}
	Which classifiable $C^*$-algebras can be expressed as crossed products of groups actions on compact metrizable spaces? 
	In other words, when are they isomorphic to groupoid $C^*$-algebras of transformation groupoids? 
\end{quesx}
It is indeed known that UHF algebras, or more generally, AF algebras, can be represented as crossed products of minimal free actions of locally finite groups on the Cantor set, as mentioned in \cite[Example 8.1.24]{GKPT}. Moreover, in \cite[Proposition 3.24]{Ma2022Purely}, it was shown that $\mathcal{O}_2$ can be expressed as a minimal purely infinite dynamical system of $\mathbb{Z}_2*\mathbb{Z}_3$. But in general, this question is still widely open. 

A prominent unresolved case regarding this question is the Jiang-Su algebra $\mathcal{Z}$.
It is well-known that $\mathcal{Z}$ cannot be realized as a $\mathbb{Z}$-crossed product using Pimsner-Voiculescu exact sequences.  In this paper, using the bijectivity of the Baum-Connes assembly for amenable groups and $KK$-theory with real coefficients, we demonstrate that there is an obstruction on the acting group for $\CZ$ to be written in the form of a crossed product. See Theorem \ref{thm: Jiang Su crossed product} below.

\begin{thmx}\label{thm: Jiang Su realizability}
Suppose $\CZ\cong C(X)\rtimes_r \Gamma$ for an action $\alpha$ of countable discrete group $\Gamma$ on a compact metrizable space $X$. Then $\Gamma$ has to be amenable, torsion free, and rationally acyclic.
\end{thmx}

Since these conditions on $\Gamma$ rules out almost all naturally occurring examples, 
this result supports the view that in the study of classifiable $C^*$-algebras, it is more convenient to work with groupoids rather than only group actions. We also refer to \cite[P. 46]{STW2025} for more discussion on Theorem \ref{thm: Jiang Su realizability}.

 \section{Preliminaries}
In this section, we recall some basic backgrounds on locally compact Hausdorff \'{e}tale groupoids and their $C^*$-algebras. We refer to \cite{Renault} and \cite{Sims} as standard references for these topics. Using the terminology in \cite{Sims}, we denote by $\GU$ the  \textit{unit space} of $\CG$. We write \emph{source} and \emph{range} maps $s,r:\CG\rightarrow\GU$, respectively and they are defined by $s(\gamma)=\gamma^{-1}\gamma$ and $r(\gamma)=\gamma\gamma^{-1}$ for $\gamma\in\CG$.
When a groupoid $\CG$ is endowed with a locally compact Hausdorff
topology under which the product and inverse maps are continuous,
the groupoid $\CG$ is called a locally compact Hausdorff groupoid.
A locally compact Hausdorff groupoid $\CG$ is called \textit{étale}
if the range map $r$ is a local homeomorphism from $\CG$ to itself,
which means for any $\gamma\in\CG$ there is an open neighborhood
$U$ of $\gamma$ such that $r(U)$ is open and $r|_{U}$ is a homeomorphism.
A set $B$ is called an \textit{$s$-section} (respectively, $r$-section) 
if there is an open set $U$ in $\CG$ such that $B\subseteq U$ and
the restriction of the source map $s|_{U}:U\rightarrow s(U)$ (respectively,
the range map $r|_{U}:U\rightarrow r(U)$) on $U$ is a homeomorphism
onto an open subset of $\GU$. The set $B$ is called a \textit{bisection} if it is both an $s$-section and a $r$-section at the same time. It is not hard to see a locally compact
Hausdorff groupoid is étale if and only if its topology has a basis
consisting of open bisections. We say a locally compact Hausdorff
étale groupoid $\CG$ is \textit{ample} if its topology has a basis
consisting of compact open bisections.

For any set $D\subseteq\GU$, Denote by 
\[
\CG_{D}\coloneqq\{\gamma\in\CG:s(\gamma)\in D\},\ \CG^{D}\coloneqq\{\gamma\in\CG:r(\gamma)\in D\},\ \text{and}\ \ \CG_{D}^{D}\coloneqq\CG^{D}\cap\CG_{D}.
\]
Note that $\CG|_{D}\coloneqq\CG_{D}^{D}$ is a subgroupoid of $\CG$
with the unit space $D$. For the singleton case $D=\{u\}$, we write $\CG_{u}$, $\CG^{u}$
and $\CG_{u}^{u}$ instead for simplicity. In this situation, we call
$\CG_{u}$ a \textit{source fiber} and $\CG^{u}$ a \textit{range
	fiber}. In addition, each $\CG_{u}^{u}$ is a group, which is called
the \textit{isotropy} at $u$. We also denote by \[\opIso(\CG)=\bigcup_{u\in \GU}\CG^u_u=\{x\in \CG: s(x)=r(x)\}\] the isotropy of the groupoid $\CG$. We say a groupoid $\CG$ is \textit{principal}
if $\opIso(\CG)=\GU$. A groupoid $\CG$ is called \textit{topologically
	principal} if the set $\{u\in\GU:\CG_{u}^{u}=\{u\}\}$ is dense in
$\GU$. The groupoid $\CG$ is also said to be \textit{effective} if $\opIso(\CG)^o=\GU$. Recall that effectiveness is equivalent to the topological principality if $\CG$ is second countable (see \cite[Lemma 4.2.3]{Sims}). A subset $D$ in $\GU$ is called $\CG$-\textit{invariant}
if $r(\CG D)=D$, which is equivalent to the condition $\CG^{D}=\CG_{D}$.
A groupoid $\CG$ is called \textit{minimal} if there are no proper
non-trivial closed $\CG$-invariant subsets in $\GU$.

We now recall a general framework   introduced in \cite[Section 3]{M-W} on castles and towers in a groupoid. 

\begin{defn} \label{defn:plurisection}
	Let $\CG$ be a groupoid. 
	A \emph{transverse system} 
	in $\CG$ is a disjoint collection $\CC$ of nonempty subsets in $\CG$ satisfying
	the following conditions: 
	\begin{enumerate}[label=(\roman*)]
		\item \label{defn:plurisection:unary} for any $C \in \CC$, we have $r(C), s(C), C^{-1} \in \CC$, and
		\item \label{defn:plurisection:binary} for any $C, D \in \CC$, we have $C D \in \CC \cup \{\varnothing\}$. 
	\end{enumerate}
	Given a transverse system $\CC$, we write $\CC^{(0)} := \left\{ C \in \CC \colon C \subseteq \CG^{(0)} \right\}$. 
\end{defn}
It is direct to see that  for a transverse system $\CC$, the union $\bigcup \CC$ is a groupoid with the unit space $\CC^{(0)}$ (see \cite[Lemma 3.2(5)]{M-W}). We refer to \cite[Lemma 3.2]{M-W} for additional basic properties for transverse systems.

\begin{defn} \label{defn:invariant-functions}
	Let $\CG$ be a groupoid and let $\CH$ be a subgroupoid of $\CG$. Define the equivalence relation $\sim_{\CH}$ on $\CG$ such that $x \sim_{\CH} y$ if and only if there are $z, w \in \CH \cup \CG^{(0)}$ such that $x = z y w$. 
\end{defn}

\begin{defn}\label{defn: towers and castles} 
	A \emph{castle} in $\CG$ is a transverse system $\CC$ that is finite and principal in the sense that the map $r \times s \colon \CC \to \CC^{(0)} \times \CC^{(0)}$, $C \mapsto (r(C) , s(C))$ is injective. 
	
	The members of $\CC^{(0)}$ are sometimes called the \emph{levels} of the castle $\CC$ or simply \emph{$\CC$-levels}, while
	the members of $\CC$ are sometimes called the \emph{ladders} of the castle $\CC$ or simply \emph{$\CC$-ladders}. 
	
	We say a castle is \emph{open} (respectively, \emph{closed}, \emph{clopen}, \emph{Borel}, \emph{compact} or \emph{precompact}) if each of its members is \emph{open} (respectively, \emph{closed}, \emph{clopen}, \emph{Borel}, \emph{compact} or \emph{precompact}) in $\CG$. 
			
	Let $\sim_{\CC}$ be the equivalence relation. Then the equivalence classes of $\sim_{\CC}$ are called the \emph{towers} in $\CC$ or simply \emph{$\CC$-towers}. 
	In particular, if $\sim_{\CC}$ has at most one equivalence class, then $\CC$ itself may also be called a \emph{tower}. 
\end{defn}

One concrete way to interpret the tower is the
 \emph{multisection}, that was introduced by Nekrashevych in \cite[Definition~3.1]{Ne}. Let $K\subseteq \CG$ be a compact set and $\CC$ a castle. We define the \emph{$K$-shrinking} of $\CC$ by 
	\begin{align*}
		\CC^{-K} := \left\{ C \in \CC \colon \left( K \cup K^{-1} \right) C \subseteq \bigcup \CC \text{ and }  C \left( K \cup K^{-1} \right) \subseteq \bigcup \CC \right\} \; .
	\end{align*}
which is  also a castle by \cite[Lemma 5.2]{M-W}.
\begin{defn}
    Let $K\subseteq \CG$ be a compact set. We say a castle $\CC$ is $K$-\textit{extendable} to a castle $\CD$ if $\CC\subseteq \CD^{-K}$.
\end{defn}

\begin{rmk}\label{rmk: extendibility}
    Suppose the compact set $K$ is symmetric, i.e., $K=K^{-1}$. Then, using the fact $\CC$ is a castle, it is straightforward to verify (see \cite[Remark 5.4]{M-W}) that $\CC$ is $K$-extendable to $\CD$ if and only if $K\cdot \bigcup \CC\subseteq \bigcup \CD$.
\end{rmk}

We then recall the groupoid subequivalence relation. 
\begin{defn}[{\cite[Section 7]{M-W}}] \label{defn: groupoid subequivalence}
    Let $\CG$ be a groupoid, $K$ be a subset in $\GU$ and $U, V$ be open subsets in $\GU$. We write
    \begin{enumerate}[label=(\roman*)]        \item $K\prec_\CG U$ if there 
    is an open $s$-section $A$ such that $K \subseteq r(A)$ and $s(A) \subseteq U$. 
        \item $U\precsim_\CG V$ if $K\prec_\CG V$ holds for every compact subset $K\subseteq U$.

        \item $U\precsim_{\CG, 2}V$ if for any compact $K\subseteq U$ there are disjoint non-empty open subsets $V_1, V_2\subseteq V$ such that $K\prec_\CG V_1$ and $K\prec_\CG V_2$. 
        \end{enumerate}
\end{defn}

\begin{rmk}
\begin{enumerate}[label=(\roman*)]    \item We remark that the relation "$K\prec_\CG U$" is usually defined in the situation that $K$ is compact. See, e.g., \cite{Ma2022Purely} and \cite{M-W}. However, it is harmless to extend the definition to a general set $K$ in $\GU$, which simplifies certain calculations in Proposition \ref{prop: almost elementary partial dyna sys} by using Lemma \ref{lem: simplify groupoid subequivalence}.
    \item We would like to point out that the notion of ``$K\prec_\CG U$'' as defined in Definition \ref{defn: groupoid subequivalence}(i) is strictly weaker than the relation carrying the same notation in \cite[Section 3]{Ma2022Purely}, where the open $s$-section $A$ is further required to be a disjoint union of open bisections. However, it is worth noting that the two notions coincide in the case of ample groupoids (by refining open covers to clopen partitions; see \cite[Remark~10.5]{M-W}) and transformation groupoids of partial dynamical systems (see Remark~\autoref{rmk:dynamical_subequivalence_relation} below). 
    Since our main examples in the following discussion fall within these two classes of groupoids, 
    this subtle distinction does not make a significant difference for the results in this paper. 
    See \cite[Remark~7.4]{M-W} for a more in-depth discussion of this distinction. 
    \end{enumerate}
\end{rmk}

We then recall the definition of almost elementariness of groupoids with compact unit spaces.

\begin{defn}\cite[Definition 7.6]{M-W}\label{defn: AE} 
	Let $\CG$ be a locally compact Hausdorff \'{e}tale groupoid with a compact unit space. 
	We say that $\CG$ is \emph{almost elementary}
	if for any compact set $K$, 
	any $\varepsilon > 0$, 
	 any finite open cover $\CV$ of $\GU$, and any non-empty open set $O$ in $\GU$,
	there is a non-empty open castle $\CC$ 
	satisfying 
	\begin{enumerate}[label=(\roman*)]
		\item $\CC$ is $K$-extendable to an open castle $\CD$; 
		\item $\CC^{(0)}$ refines $\CV$; 
		\item   $\CG^{(0)}\setminus( \bigcup \CC^{(0)})\prec_{\CG}O$. 
	\end{enumerate}
\end{defn}

Finally, throughout the paper, we write $B\sqcup C$ to indicate that
the union of sets $B$ and $C$ is a disjoint union. In addition,
we denote by $\bigsqcup_{i\in I}B_{i}$ for the disjoint union of
the family $\{B_{i}:i\in I\}$. From now on, we only consider \textbf{locally compact, Hausdorff, $\sigma$-compact, étale, topological groupoids.  The word ``groupoid'' below is reserved for this kind of groupoids.} We remark that we do NOT always assume the unit space $\GU$ is compact.

\section{Almost elementary groupoids of partial dynamical systems}\label{sec: AE model for partial ds}
In this section, we study transformation groupoids of partial dynamical systems, which form groupoid models for stably finite $C^*$-algebras. We first recall the definition of partial dynamical systems and refer to \cite{Exel} as a standard reference for this topic.

\begin{defn}\cite[Definition 2.1, Proposition 2.5]{Exel}
    A \textit{partial action} $\beta$ of a group $\Gamma$ on a set $X$  is a pair $\beta=(\{D_g\}_{g\in \Gamma}, \{\beta(g)\}_{g\in \Gamma})$, in which $\beta(g)$ are maps $D_{g^{-1}}\to D_{g}$ such that 
    \begin{enumerate}[label=(\roman*)]        \item $D_{e_\Gamma}=X$ and $\beta(e_\Gamma)$ is the identity map;
        \item $\beta(g)(D_{g^{-1}}\cap D_h)\subseteq D_{gh}$ for any $g, h\in \Gamma$;
        \item $\beta(g)(\beta(h)(x))=\beta(gh)(x)$ for any $g, h\in \Gamma$ and $x\in D_{h^{-1}}\cap D_{(gh)^{-1}}$.
    \end{enumerate}
\end{defn}
In the topological dynamical setting throughout this section, we additionally assume the following throughout the section.
\begin{enumerate}[label=(\roman*)] 
\item The acting group $\Gamma$ is countable and discrete.
\item The underlying space $X$ is a compact metrizable space.
    \item Each $D_g$ is an open set in $X$.
    \item Each $\beta(g): D_{g^{-1}}\to D_g$ is a homeomorphism. 
    \end{enumerate}

One can define groupoids from partial dynamical systems. Similar to the global action case, we write $gx$ for $\beta(g)(x)$ for simplicity if the context is clear.

\begin{Exl}
	\label{exa:transformation-groupoid}Let $X$ be a compact
	Hausdorff space and $\Gamma$ be a discrete group. Then any partial action
	$\beta:\Gamma\curvearrowright X$ by (partial) homeomorphisms induces a locally compact
	Hausdorff étale groupoid 
	\[
	X\rtimes_\beta\Gamma\coloneqq\{(g x,g,x):g\in\Gamma,x\in D_{g^{-1}}
 \}
	\]
	equipped with the relative topology as a subset of $X\times\Gamma\times X$.
	In addition, $(g x,g,x)$ and $(h y,h,y)$ are composable
	only if $h y=x$ and 
	\[
	(g x,g,x)(h y,h,y)=(ghy,gh, y).
	\]
	One also defines $(g x,g,x)^{-1}=(x,g^{-1},g x)$
	and declares that $\GU\coloneqq\{(x,e_{\Gamma},x):x\in X\}$. It
	is not hard to verify that $s(g x,g,x)=x$ and $r(g x,g,x)=g x$.
	The groupoid $X\rtimes_\beta\Gamma$ is called a \textit{partial transformation
		groupoid}. 
\end{Exl}

\begin{rmk}\label{rmk:dynamical_subequivalence_relation}
 Let $\beta: \Gamma\curvearrowright X$ be a partial transformation groupoid $\CG=X\rtimes_\beta \Gamma$ always comes with a clopen partition 
    \[
	X\rtimes\Gamma = \bigsqcup_{g\in\Gamma} \{(g x,g,x): x\in D_{g^{-1}}\} \; ,
    \]
    thanks to the discreteness of $\Gamma$. 
    Combining this with a compactness argument, similar to \cite[Example 7.3]{M-W}, it is straightforward to see that for any subsets $K, U$ in $\GU$ with $U$ open, we have $K\prec_\CG O$ if and only if there are open sets $U_1,\dots, U_n$ in $X$ and group elements $g_1\dots, g_n\in \Gamma$ such that 
\begin{enumerate}[label=(\roman*)]     \item $U_i\subseteq D_{g^{-1}_i}$ for all $i=1,\dots, n$;
    \item $K\subseteq \bigcup_{i=1}^nU_i$ and $\bigsqcup_{i=1}^ng_iU_i\subseteq O$.
    \end{enumerate}    
    For the ease of our presentation, we often write $K\prec_\beta O$ instead of $K\prec_\CG O$. 
    This is a straightforward generalization of the notion of dynamical subequivalence for global actions (see, e.g., \cite{D}) to the setting of partial actions. 
\end{rmk}

The following lemma is straightforward but useful.

\begin{lem}\label{lem: simplify groupoid subequivalence}
    Suppose $\beta$ is a partial action of $\Gamma$ on $X$. Let $K_1, \dots, K_n$ be sets in $X$ and  $O_1,\dots, O_n$ non-empty disjoint open sets in $X$. Suppose $K_i\prec_\beta O_i$ holds for any $i=1,\dots, n$. Then $\bigcup_{i=1}^nK_i\prec_\beta \bigsqcup_{i=1}^nO_i$.
\end{lem}

We then describe towers and castles in the setting of partial dynamical systems.

\begin{defn}\label{defn: tower for partial action}
Let $\beta: \Gamma\curvearrowright X$ be a partial action of $\Gamma$ on a compact Hausdorff space. Let $S$ be a finite subset of $\Gamma$ and $B$ an open set in $X$. We say $(S, B)$ is an \textit{open tower} if 
\begin{enumerate}[label=(\roman*)]    \item $B\subseteq \bigcap_{s\in F}D_{s^{-1}}$;
    \item $\beta(s)(B)\subseteq D_{st^{-1}}$ for any $s, t\in S$
    \item $\{\beta(s)(B):s\in S\}$ is a disjoint family.
\end{enumerate}
A finite collection $\{(S_i, B_i): i\in I\}$ of towers is called an \textit{open castle} if $S_iB_i\cap S_jB_j=\emptyset$ for any distinct $i, j\in I$.
\end{defn}

In the following sense, the towers and castles for a partial dynamical system coincide with the same notions in the transformation groupoid (Definition \ref{defn: towers and castles}).

\begin{rmk}\label{rem: different notions of castles are same}
    Let $\CG=X\rtimes_\beta \Gamma$ be the partial transformation groupoid of a partial action $\beta: \Gamma\curvearrowright X$ and $(S, B)$ a tower in $\CG$. Then the collection
\[\CT(S, B)=\{(tB\times\{ts^{-1}\}\times sB)\cap \CG: s, t\in S\}\]
is an open tower in the sense of Definition \ref{defn: towers and castles}. Now if $C=\{(S_i, B_i): i\in I\}$ is an open castle in the sense of Definition \ref{defn: tower for partial action}, then the  corresponding family $\CC(C)=\{\CT(S_i, B_i): i\in I\}$ is an open castle in $\CG$ in the sense of Definition \ref{defn: towers and castles}. 
\end{rmk}

\begin{lem}\label{lem: extendable-partial action}
Let $\CG=X\rtimes_\beta \Gamma$ be the transformation groupoid of a partial action $\beta: \Gamma\curvearrowright X$. Let $C =\{(S_i, B_i): i\in I\}$ and $C' =\{(T_i, B_i): i\in I\}$ be two open castles with the same base $B_i$ and with $S_i\subseteq T_i$ for each $i\in I$. Denote by $\CC=\CC(C)$ and $\CD=\CC(C')$ the associated open castles defined above. Let $K$ be a symmetric finite subset of $\Gamma$, which induces a set 
\[M=\{(gx, g, x): g\in K, x
\in D_{g^{-1}}\}\]
in $\CG$. Suppose $KS_i\subseteq T_i$ for each $i\in I$. Then $\CC$ is $M$-extendable to $\CD$.
\end{lem}
\begin{proof}
Since $K$ is symmetric, so is $M$. Let $z_1=(tx, ts^{-1}, sx)\in \CT(S_i, B_i)$ for some $i\in I$ and $z_2=(gy, g, y)\in M$. Suppose $z_2\cdot z_1$ is composable in $\CG$, which means $y=tx\in D_{g^{-1}}$. Then one has $z_2\cdot z_1=(gtx, gts^{-1}, sx)$. Now, $gt\in KS_i\subseteq T_i$ implies that $z_2\cdot z_1=(t'x, t's^{-1}, sx)$ for some $t'=gt\in T_i$ and thus $z_2\cdot z_1\in \CT(T_i, B_i)$. This shows that $\CC$ is $M$-extendable to $\CD$ by Remark \ref{rmk: extendibility}.
\end{proof}

A natural way to define a partial dynamical system is to restrict a (global) group action to a certain family of open sets (see, e.g., \cite[section 3]{Exel}).  Let $\alpha: \Gamma\curvearrowright X$ be a global action of a countable amenable discrete group $\Gamma$ on a compact Hausdorff space $X$. We now focus on minimal partial dynamical systems $\beta: \Gamma\curvearrowright X$ obtained by restricting a global minimal action $\alpha$ to certain open sets $D_{g^{-1}}=X\setminus \bigcup_{j\in I_g} Y_{j, g}$, where $J_g$ is a finite index set (with $J_e = \emptyset$) and each $Y_{j, g}$ is compact set in $X$ intersecting every $\alpha$-orbit at most once. We first record some basic properties of such sets.

\begin{rmk}\label{rem: basic property of topo small set}
  Let $\alpha: \Gamma\curvearrowright X$ be a minimal free global action by an infinite group. Let $C$ be a closed set in $X$ that intersects every $\alpha$-orbit at most once. Thus, for any $s\in \Gamma$, the translation $sC$ intersects any $\alpha$-orbit at most once as well. In addition, it is direct to see that $\sup_{\mu\in M_{\alpha, \Gamma}(X)}\mu(C)=0$ because $\{sC: s\in \Gamma\}$ is a disjoint family. In particular, any non-empty open set cannot be covered by countably many such compact sets intersecting any $\alpha$-orbit at most once.
  \end{rmk}

We also need the following lemma.

\begin{lem}\label{lem: topo small set subsequivalent}
Let $\alpha: \Gamma \curvearrowright X$ be a global minimal free action. Let $C$ be a compact set that meets each $\alpha$-orbit at most once and $O$ a non-empty open set in $X$. Suppose $\beta: \Gamma \curvearrowright X$ is a minimal partial action obtained by a restriction of $\alpha$. Then one has $C\prec_\beta O$ in the transformation groupoid $X\rtimes_\beta \Gamma$. 
\end{lem}
\begin{proof}
    Since $C$ meets every $\alpha$-orbit at most once, one has $\{\alpha(g)(C): g\in \Gamma\}$ is a disjoint family. Then the minimality of $\beta$ and compactness of $X$ imply that there are a finite family $F\subseteq \Gamma$ and open sets $O_g$ for $g\in F$ such that 
    \begin{enumerate}[label=(\roman*)]        \item $O_g\subseteq D_{g^{-1}}$ for any $g\in F$;
        \item $\{O_g: g\in F\}$ forms an open cover of $X$;
        \item $\bigcup_{g\in F}\beta(g)(O_g)\subseteq O$.
    \end{enumerate}
   Because $X$ is  Hausdorff, one chooses an open set $U\supset C$ such that $\{\alpha(g)(U): g\in F\}$ is a disjoint family. Now observe that 
    \begin{enumerate}[label=(\roman*)]        \item $C\subseteq \bigcup_{g\in F}O_g\cap U$,
        \item $O_g\cap U\subseteq D_{g^{-1}}$, and
        \item $\bigsqcup_{g\in F}\beta(g)(O_g\cap U)\subseteq O$.
    \end{enumerate}
    This implies that $C\prec_\beta O$.
\end{proof}

Using the almost finiteness introduced in \cite[Definition 8.2]{D} in the setting of global dynamical systems, we have the following. Recall that for a minimal partial action  $\beta: \Gamma\curvearrowright X$ of an infinite countable group $\Gamma$ on an infinite compact Hausdorff space $X$, the space $X$ is necessarily \textit{perfect} in the sense that  $X$ has no isolated points. Suppose not and let $\{x\}$ be open in $X$. Then the minimality of $\beta$ implies that for any $x'\in X$ there exists a $g\in \Gamma$ such that $\beta(g)(x')\in \{x\}$. This implies that the orbit $\{\beta(g)(x): x\in D_{g^{-1}}, g\in \Gamma\}$ covers $X$ and thus $X$ is finite by compactness. But this is a contradiction to the setting.

\begin{prop}\label{prop: almost elementary partial dyna sys}
Let $\alpha: \Gamma\curvearrowright X$ be a minimal free global action, which is almost finite in the sense of  \cite[Definition 8.2]{D}. Suppose $\beta:\Gamma\curvearrowright X$ is a minimal partial dynamical system obtained by restricting $\alpha$ to a family $\{D_g: g\in \Gamma\}$ of open sets in $X$ such that each $D_{g^{-1}}$ is of the form $D_{g^{-1}}=X\setminus \bigcup_{j\in J_g}Y_{j, g}$, where $J_g$ is a finite index set and each $Y_{j, g}$ is a compact set meeting any $\alpha$-orbit at most once. Then the groupoid $X\rtimes_\beta \Gamma$ is almost elementary.
\end{prop}
\begin{proof}
    Denote by $\CG=X\rtimes_\beta \Gamma$ the groupoid of the partial action $\beta$. Let $C$ be a compact set in $\CG$, $O$ a non-empty open set in $X=\GU$, and  $\CV$ an open cover of $X$. First choose finite set $K\subseteq \Gamma$ such that 
\[C\subseteq \{(\beta(g)(x), g, x): g\in K, x\in D_{g^{-1}}\}.\]
Since $\alpha$ is minimal, one has $\varepsilon=\inf_{\mu\in M_{\Gamma, \alpha}(X)}\mu(O)>0$. Using the almost finiteness of $\alpha$, it is direct to see that there exists a finite open castle $\{(T_i, V_i): i\in I\}$ for the global action $\alpha$ such that 
\begin{enumerate}[label=(\roman*)]    \item all $T_i$ are $(K, \varepsilon/3)$-invariant in the sense that $|\bigcap_{t\in K}t^{-1}T_i|\geq (1-\varepsilon/3)|T_i|$;
    \item for any $t\in T_i$, one has $\alpha(t)(V_i)\subseteq V$ for some $V\in \CV$;
    \item $\sup_{\mu\in M_{\Gamma, \alpha}(X)}\mu(X\setminus \bigsqcup_{i\in I, t\in T_i}\alpha(t)(V_i))<\varepsilon/3$.
    \end{enumerate}
Without loss of any generality, we may assume $e_\Gamma\in T_i$ for each $i\in I$ by shifting $T_i$ on the right if necessary. Now, write $S_i=\bigcap_{t\in K}t^{-1}T_i$ and $R=X\setminus \bigsqcup_{i\in I, t\in S_i}\alpha(t)(V_i)$ for simplicity. Observe that $\sup_{\mu\in M_{\Gamma, \alpha}(X)}\mu(R)<2\varepsilon/3$. On the other hand, the almost finiteness for $\alpha$ implies that $\alpha$ has dynamical strict comparison in the sense of \cite[Definition 3.2]{D} by \cite[Theorem 9.2]{D}, i.e., there is a finite set $F\subseteq \Gamma$ and a family $\{O_g: g\in F\}$ of open sets in $X$ such that $R\subseteq \bigcup_{g\in F}O_g$ and $\bigsqcup_{g\in F}\alpha(g)(\overline{O_g})\subsetneq O$ by the normality of the space $X$. Write $O'=\bigsqcup_{g\in F}\alpha(g)(O_g)$ for simplicity.

Now, denote by $\CK_g$ the family of compact sets $Y_{j, g}\subseteq D^c_{g^{-1}}$ satisfying $Y_{j, g}\cap O_g\neq \emptyset$. Also, for any $i\in I$ and $s\in T_i$, write $\CL_{i, s}$ the collection of all compact sets $Y_{j, ts^{-1}}\subseteq D^c_{st^{-1}}$ satisfying $Y_{j, ts^{-1}}\cap \alpha(s)(V_i)\neq \emptyset$ for some $t\in T_i$.

Since $X$ is Hausdorff and perfect, one chooses $k=\sum_{i\in I, s\in T_i}|\CL_{i, s}|\cdot|T_i|+\sum_{g\in F}|\CK_g|$ disjoint open sets  contained in $O\setminus \bigsqcup_{g\in F}\alpha(g)(\overline{O_g})$ and we denote $\CU$ for this family of open sets for simplicity.  
For every $g\in F$, choose a subfamily $\CU_g\subseteq \CU$ with cardinality $|\CK_g|$ and satisfying $\CU_g\cap \CU_h=\emptyset$ if $h\neq g\in F$. Fix
an arbitrary bijective map $\varphi_g$ from $\CK_g$ to $\CU_g$.
Since the groupoid $X\rtimes_\beta \Gamma$ is minimal, Lemma \ref{lem: topo small set subsequivalent} implies that $Y_{j, g}\prec_\beta \varphi_g(Y_{j, g})$ for any $Y_{j, g}\in \CK_g$.   On the other hand, choose subfamilies $\CV_{i, s}\subseteq \CU$ for $s\in T_i$ such that 
\begin{enumerate}[label=(\roman*)] \item $|\CV_{i, s}|=|\CL_{i, s}|\cdot |T_i|$ for any $i\in I, s\in T_i$;
    \item $\CV_{i, s}\cap \CV_{l, t}=\emptyset$ if $(i, s)\neq (l, t)$;
    \item $\CV_{i, s}\cap \CU_g=\emptyset$ for any $i\in I, s\in T_i, g\in F$.
\end{enumerate}
This is possible because $|\CU|=k$. Now, for each $i\in I, s\in T_i$ and $Y\in \CL_{i,s}$, the translation $\alpha(ts^{-1})(Y)$ still intersects any $\alpha$-orbit at most once for any $t\in T_i$. Fix an arbitrary bijective map $\psi_{i, s}: \CL_{i, s}\times T_i\to \CV_{i, s}$ and one has $\alpha(ts^{-1})(Y)\prec_\beta \psi_{i,s}(Y, t)$ by Lemma \ref{lem: topo small set subsequivalent}.

Then, for each $i\in I$, define $W_i=V_i\setminus \bigcup\{\alpha(s^{-1})(Y): Y\in \CL_{i, s}, s\in T_i\}$, which is a non-empty open set by Remark \ref{rem: basic property of topo small set}. Observe that $\CC=\{(S_i, W_i): i\in I\}$ and $\CD=\{(T_i, W_i): i\in I\}$ are $\beta$-castles in the sense of Definition \ref{defn: tower for partial action} by the construction. Abusing notations a bit, we also denote by $\CC$ and $\CD$ the groupoid they induced by Remark \ref{rem: different notions of castles are same}. Now because $KS_i\subseteq T_i$ holds for each $i\in I$, Lemma \ref{lem: extendable-partial action} implies that $\CC$ is $C$-extendable to $\CD$. In addition, one still has that for any $t\in T_i$, the $\CD$-level $tW_i\subseteq V$ for some $V\in \CV$.

Then, decompose $R= (R\cap \bigcap_{g\in F}D_{g^{-1}})\cup (\bigcup_{g\in F}R\cap D^c_{g^{-1}})$.
First, one has
\[(R\cap \bigcap_{g\in F}D_{g^{-1}})\subseteq \bigcup_{g\in F}(O_g\cap D_{g^{-1}})\]
while
\[\bigsqcup_{g\in F}\beta(g)(O_g\cap D_{g^{-1}})\subseteq O'.\]
This means $R\cap \bigcap_{g\in F}D_{g^{-1}})\prec_\beta O'$.
In addition, observe 
\[\bigcup_{g\in F}(R\cap D^c_{g^{-1}})\subseteq \bigcup_{g\in F}(\bigcup \CK_g).\]
Recall that for any $g\in F$ and $Y_{j, g}\in \CK_g$, one has $Y_{j, g}\prec_\beta \varphi_g(Y_{j,g})\in \CU_g$. Lemma \ref{lem: simplify groupoid subequivalence} then implies that 
$R\prec_\beta O'\sqcup \bigsqcup_{g\in F}\bigsqcup\CU_g$. Finally, note that
\[X\setminus \bigsqcup_{i\in I}\bigsqcup_{s\in S_i}\beta(s)(W_i)\subseteq R\cup \bigcup\{\alpha(ts^{-1})(Y): Y\in \CL_{i, s}, s\in T_i, i\in I \},\]
which actually implies that 
\[X\setminus \bigsqcup_{i\in I}\bigsqcup_{s\in S_i}\beta(s)(W_i)\prec_\beta O'\sqcup \bigsqcup_{g\in F}\bigsqcup\CU_g\sqcup \bigsqcup_{i\in I, s\in T_i}\bigsqcup \CV_{i, s}\subseteq O\]
by Lemma \ref{lem: simplify groupoid subequivalence}. Thus, one has $\CG$ is almost elementary.
\end{proof}

The following is a consequence of  Proposition \ref{prop: almost elementary partial dyna sys} by using \cite[Corollary~13.13]{M-W}. 

\begin{thm}\label{thm: almost elementary partial action model}
    Let $\beta:\Gamma\curvearrowright X$ be a partial action in Proposition \ref{prop: almost elementary partial dyna sys}. Denote by $\CG=X\rtimes_\beta \Gamma$. Then $\CG$ is a minimal, amenable, almost elementary, second countable groupoid on a compact space. Therefore, the reduced groupoid $C^*$-algebra $C^*_r(\CG)$ is unital simple nuclear separable $\CZ$-stable and thus classifiable by its Elliott invariant.
\end{thm}

Now, we turn to the case $\alpha: \Z\curvearrowright X$ be a global minimal (free) action induced by a homeomorphism $\varphi: X\to X$. For each integer $n$, write 
\[
    I_n \coloneqq ((-\infty, 0] \Delta (-\infty, n]) \cap \Z = 
    \begin{cases}
        [1, n]\cap \Z \, , & n > 0 \\
        \emptyset \, , & n = 0 \\
        [n+1, 0]\cap \Z \, , & n < 0 
    \end{cases}
    \; .
\]
A simple but useful observation for any integers and $n, m\in \Z$ are following.
   \begin{enumerate}
       \item $I_{-n}+n=I_n$ and
       \item $I_{n+m}-n\subseteq I_{-n}\cup I_m$.
   \end{enumerate}
Let $Y$ be a closed set in $X$ meeting each $\alpha$-orbit at most once. For any $g\in \Z$, define $D_g=X\setminus \bigsqcup_{n\in I_g}\varphi^n(Y)$. Then define partial action $\beta$ by restricting $\alpha$ on all of these $D_{g^{-1}}$. One may verify directly that the partial system $\beta$ is well-defined. Since $X$ is compact, the one-sided half $\alpha$-orbits $\{\varphi^n(x): n\geq 1\}$ and $\{\varphi^n(x): n \leq 0\}$ are both dense in $X$ for any $x\in X$, from which one concludes that the partial action $\beta$ is minimal as well. In addition, $\beta$ is free because $\alpha$ is free. Then the transformation groupoid $\CG=X\rtimes_\beta\Z$ is minimal and principal. Actually, the transformation groupoid $\CG=X\rtimes_\beta\Z$ for $\beta$ is nothing but the equivalence relation
\[\CR_Y=\{(x, \varphi^n(x)): n\in \Z, x\in X\}\setminus \{(\varphi^k(y), \varphi^l(y)): y\in Y, l<1\leq k\ \text{or }k<1\leq l\},\]
which is exactly the so-called \textit{orbit-breaking} equivalence relation considered in \cite{DPS2}. For this kind of groupoids, we have the following result as a corollary of Proposition \ref{prop: almost elementary partial dyna sys} and Theorem \ref{thm: almost elementary partial action model}.

\begin{cor}\label{cor: almost elementary partial dyn}
    Let $\alpha: \Z\curvearrowright X$ be a minimal global action induced by a homomorphism $\varphi$ on a compact metrizable space $X$ with the finite covering dimension. Let $Y$ be a compact set in $X$ intersecting any $\alpha$-orbit at most once and $I_g$ be intervals of integers defined above for all $g\in \Z$. Suppose $\beta: \Z\curvearrowright X$ is the partial action obtained by restricting $\alpha$ to the family $\{D_g=X\setminus \bigsqcup_{n\in I_g}\varphi^n(Y): g\in \Z\}$. Then the transformation groupoid $\CG=X\rtimes_\beta \Gamma$ is minimal, amenable, almost elementary, second countable groupoid and thus the groupoid $C^*$-algebra $C^*_r(\CG)$ is $\CZ$-stable and thus classifiable by its Elliott invariant.
\end{cor}
\begin{proof}
    In light of Proposition \ref{prop: almost elementary partial dyna sys} and Theorem \ref{thm: almost elementary partial action model}, it suffices to show $\alpha$ is almost finite in the sense of \cite[Definition 8.2]{D}. But this has been established in \cite[Theorem C]{K-S}.
\end{proof}

It was proved in \cite{DPS1} and \cite{DPS2} that the orbit-breaking equivalence relation $R_Y$ defined above provides groupoid model for a large class of classifiable $C^*$-algebras including Jiang-Su algebra $\CZ$. 

\begin{thm}\cite[Corollary 6.4]{DPS2}\label{thm: DPS groupoid model}
    Let $G_0$ and $G_1$ be countable abelian groups, $k\in \Z_{>0}$, $\Delta$ a finite-dimensional Choquet simplex. Then the pair $(\Z\oplus G_0, \Z_{>0}\oplus G_0, [1]=(k, 0))$ is an ordered abelian group, and if there is a map
    $r: \Delta\to S(\Z\oplus G_0)$ defined by $\tau\mapsto ((n, g)\mapsto n/k)$, then there exists an amenable minimal equivalence relation $R$ such that $C^*(R)$ is classifiable and
    \[\operatorname{Ell}(C^*(R))\cong ((\Z\oplus G_0, \Z_{>0}\oplus G_0, [1]=(k, 0)), G_1, \Delta, r).\]
\end{thm}

\begin{rmk}\label{rem: almost elementary amplification}
    We remark that the equivalence relation $R$ in the proof of Theorem \ref{thm: DPS groupoid model} in \cite{DPS2} is an amplification of an orbit-breaking equivalence relation $R_Y$ in the sense that $R=R_Y\times R_k$ on $X\times \{1,\dots, k\}$, where $R_k$ is the full equivalence relation on $k$ points. It is straightforward to see $R$ is almost elementary because $R_Y$ as a groupoid is almost elementary and $R_k$ is a pair groupoid. 
    \end{rmk}

Therefore, we have arrived at the main result in this section by using Corollary \ref{cor: almost elementary partial dyn}, Theorem \ref{thm: DPS groupoid model} and Remark \ref{rem: almost elementary amplification}.

\begin{thm}\label{thm: groupoid model stably finite}
   Let $G_0$ and $G_1$ be countable abelian groups, $k\in \Z_{>0}$, and $\Delta$ a finite-dimensional Choquet simplex.  Let $A$ be a  $C^*$-algebra with the Elliott invariant
\[\operatorname{Ell}(A)\cong ((\Z\oplus G_0, \Z_{>0}\oplus G_0, [1]=(k, 0)), G_1, \Delta, r),\]
in which $r: \Delta\to S(\Z\oplus G_0)$ defined by $\tau\mapsto ((n, g)\mapsto n/k)$. Then $A$ is classifiable if and only if it has a minimal, amenable, second countable, almost elementary groupoid model. In particular, Jiang-Su algebra $\CZ$ has an almost elementary groupoid model.
\end{thm}

\section{Almost elementary groupoid models for crossed product $C^*$-algebras}
In this section, as an application of the methods in section \ref{sec: AE model for partial ds}, we show that if a crossed product $A=C(Y)\rtimes_r \Gamma$, arising from a minimal free action of an amenable group $G$ on an infinite compact metrizable space $Y$ with $\dim(Y)<\infty$, is $\CZ$-stable, then $A$ has a minimal topological principal almost elementary groupoid model, which is not necessarily $Y\rtimes \Gamma$ itself.  This will complete the proof of Theorem \ref{thm: main}.

Denote by $\alpha: \Z\curvearrowright X$ the minimal free partial dynamical system serving as the groupoid model of the Jiang-Su algebra $\CZ$ in \cite{DPS1}. Recall that $X$ is of finite covering dimension. Denote by $\bar{\alpha}$ the globalization of $\alpha$, which is a minimal free action. Therefore, $\bar{\alpha}$ is almost finite by taking the usual Kakutani-Rokhlin partition. 
 Let $\beta: \Gamma\curvearrowright Y$ be a minimal action of an amenable group $G$ on an infinite compact metrizable space $Y$. Recall that in this case,  $X, Y$ are necessarily perfect.

Denote by $M_{\alpha\times \beta}(X\times Y)$ the set of all $\Z\times \Gamma$-invariant probability measures for the partial action $\alpha\times \beta$ and denote by $M_{\bar{\alpha}\times \beta}(X\times Y)$ the set of all $\Z\times \Gamma$-invariant probability measures for the global action  $\bar\alpha\times \beta$. We begin with the following lemma.
    \begin{lem}
       $M_{\alpha\times\beta}(X\times Y)=M_{\bar{\alpha}\times \beta}(X\times Y)$. 
    \end{lem}
\begin{proof}
  The direction  ``$\supset$'' is trivial. For the other direction, let $\mu\in M_{\alpha\times \beta}(X\times Y)$. First, observe $\mu(\{x\}\times Y)=0$. Indeed, there is an infinite set $N\subseteq \Z$ such that  $x\in D_{-n}$ for all $n\in N$. Note that all $\{\alpha(n)(x)\}\times Y$ for $n\in N$ form an infinite disjoint family and thus $\mu(\{x\}\times Y)=0$.

  Let $U\times V$ be an open set in $X\times Y$ and $(n, g)\in \Z\times \Gamma$. Write $D^c_{-n}=\{x_1, \dots, x_k\}$ explicitly, which is finite. Observe
  \[U\times V=((U\cap D_{-n})\times V)\sqcup\bigsqcup_{i=1}^k(\{x_i\}\times V),\]
  which implies that $\mu(U\times V)=\mu((U\cap D_{-n})\times V)$. Since $\mu\in  M_{\alpha\times \beta}(X\times Y)$, one has
  \[\mu((U\cap D_{-n})\times V)=\mu(\alpha(n)(U\cap D_{-n})\times \beta(g)(V))=\mu((\bar{\alpha}(n)(U)\cap D_n)\times\beta(g)(V)).\]
In the same way one also has $\mu((\bar{\alpha}(n)(U)\cap D_n)\times\beta(g)(V))=\mu(\bar{\alpha}(n)(U)\times\beta(g)(V))$. Thus, we have 
\[\mu(\bar{\alpha}(n)(U)\times\beta(g)(V))=\mu(U\times V),\]
which implies that $\mu\in M_{\bar{\alpha}\times\beta}(X\times Y)$.
\end{proof}

Therefore, we may write $M_{\Z\times \Gamma}(X\times Y)$ for this set of invariant probability measures without any ambiguity.

    \begin{prop}\label{prop: comparison}
        Suppose $\bar{\alpha}\times \beta: \Z\times \Gamma\curvearrowright X\times Y$ has the dynamical strict comparison, i.e., if $K$ is a compact set and $O$ is an open set satisfying $\mu(K)<\mu(O)$ for any $\mu\in M_{\Z\times \Gamma}(X\times Y)$, then $K\prec_\alpha O$ in the sense of Remark \ref{rmk:dynamical_subequivalence_relation}. Then so is the partial action $\alpha\times \beta$.
        \end{prop}
        \begin{proof}
           Let $K$ be a compact set and $O$ a non-empty open set in $X\times Y$ satisfying $\mu(K)<\mu(O)$ for all $\mu\in M_{\Z\times \Gamma}(X\times Y)$. Then by the assumption and the normality of $X\times Y$, there exist open sets $U_i\times V_i$ and group elements $(n_i, g_i)\in \Z\times \Gamma$ for each $i\in I$ such that 
           $K\subseteq \bigcup_{i\in I}(U_i\times V_i)$ and \[\bigsqcup_{i\in I}(\bar{\alpha}(n_i)(\overline{U_i})\times \beta(g_i)(\overline{V_i}))\subsetneq O.\] 
           Write $D^c_{-n_i}=\{x_{j, i}: j\in J_i\}$ for each $i\in I$, in which each $J_i$ is finite. Define $L=\sum_{i\in I}|J_i|$ and choose $L$ disjoint open sets $M_1\times N_1, \dots, M_L\times N_L$ contained in the open set  $O\setminus \bigsqcup_{i\in I}(\bar{\alpha}(n_i)(\overline{U_i})\times \beta(g_i)(\overline{V_i}))$ because $X$ and $Y$ are perfect. In addition, choose an injective map $\varphi: \{x_{j, i}: j\in J_i, i\in I\}\to \{1, \dots, L\}$. Then for $x_{j, i}$, by the minimality of $\beta$ and the compactness of $\overline{V_i}$, there exist open sets $Q_{j,i, 1}, \dots, Q_{j,i, k_{j, i}}$ in $Y$ forming a cover of $\overline{V_i}$ and $g_{j,i, 1},\dots, g_{j,i, k_{i, j}}\in \Gamma$ such that $\bigcup_{l=1}^{k_{i, j}}\beta(g_{j,i, l})(Q_{j,i, l})\subseteq N_{\varphi(x_{j, i})}$. Then the perfectness of $X$ allows to choose $k_{j, i}$ disjoint open sets $P_{j, i, 1},\dots P_{j, i, k_{j, i}}$ in $M_{\varphi(x_{j, i})}$. Now the minimality of the partial action $\alpha$ entails that there are $m_{j, i, 1}, \dots, m_{j, i, k_{j, i}}\in \Z$ and open neighborhoods $W_{j, i, 1},\dots, W_{j, i, k_{j, i}}$ of $x_{j, i}$ such that $\alpha(m_{j, i, l})(W_{j, i, l})\subseteq P_{j, i, l}$ for all $1\leq l\leq k_{j, i}$. These imply the collection $\{W_{j, i, l}\times Q_{j, i, l}: 1\leq l\leq k_{j, i}\}$ form a cover for $\{x_{j,i}\}\times V_i$ and $\{\alpha(m_{j, i, l})(W_{j, i, l})\times \beta(g_{j, i, l})(Q_{j, i, l}): 1\leq l\leq k_{j, i}\}$ is a disjoint family of open sets contained in $M_{\varphi(x_{j, i})}\times N_{\varphi(x_{j, i})}$. Therefore, the collection
           \[\{(U_i\cap D_{-n_i})\times V_i: i\in I\}\cup \{W_{j, i, l}\times Q_{j, i, l}: i\in I, j\in J_i, 1\leq l\leq k_{j, i}\}\]
           forms an open cover of $K$ and the collection
           \[\{\alpha(n_i)(U_i\cap D_{-n_i})\times \beta(g_i)(V_i),\alpha(m_{j, i, l})(W_{j, i, l})\times \beta(g_{j, i, l})(Q_{j, i, l}): i\in I, j\in J_i, 1\leq l\leq k_{j, i}\}\]
           is a disjoint family of open sets contained in $O$.
Therefore, $K\prec_{\alpha\times\beta} O$ holds.
  \end{proof}

\begin{prop}\label{prop: AE permanance}
    Suppose $\bar{\alpha}\times \beta: \Z\times \Gamma\curvearrowright X\times Y$ is almost finite. Then the partial action $\alpha\times \beta$ is almost elementary as a partial transformation groupoid.
    \end{prop}
\begin{proof}
 Denote by $\CG=(X\times Y)\rtimes_{\alpha\times \beta} (\Z\times \Gamma)$ the partial transformation groupoid of the partial action $\alpha\times \beta$. Let $K$ be a compact set in $\CG$, $O$ a non-empty open set in $X=\GU$, and  $\CV$ an open cover of $X$. First choose finite set $F\subseteq \Z\times \Gamma$ such that 
\[K\subseteq \{((\alpha(n)(x), \beta(g)(y)), (n,g), (x, y)): (n, g)\in F, x\in D_{-n}, y\in Y\}.\]
Since $\bar{\alpha}\times\beta$ is minimal, one has $\varepsilon:=\inf_{\mu\in M_{\Z\times \Gamma}(X\times Y)}\mu(O)>0$. We denote by  $\gamma=\bar{\alpha}\times \beta$ for simplicity. Then the almost finiteness of $\gamma$ implies that there exists a finite open castle $\{(T_i, V_i): i\in I\}$ for the global action $\gamma$ such that 
\begin{enumerate}[label=(\roman*)]    \item all $T_i$ are $(F, \varepsilon/3)$-invariant in the sense that $|\bigcap_{t\in F}t^{-1}T_i|\geq (1-\varepsilon/3)|T_i|$;
    \item for any $t\in T_i$, one has $\gamma(t)(V_i)\subseteq V$ for some $V\in \CV$;
    \item $\sup_{\mu\in M_{\Z\times \Gamma}(X)}\mu((X\times Y)\setminus \bigsqcup_{i\in I, t\in T_i}\gamma(t)(V_i))<\varepsilon/3$.
    \end{enumerate}
Without loss of any generality, we may still assume $e_\Gamma\in T_i$ for each $i\in I$. Now, write $S_i=\bigcap_{t\in F}t^{-1}T_i$ and $R=(X\times Y)\setminus \bigsqcup_{i\in I, t\in S_i}\gamma(t)(V_i)$ for simplicity. Observe that $\sup_{\mu\in M_{\Z\times \Gamma}(X\times Y)}\mu(R)<2\varepsilon/3$. 

Let $i\in I$ and $s, t\in T_i$. Writing $st^{-1}=(n, g)$, define
\[P_{s, t}=\{x\in D^c_{-n}: \text{there exists a }y\in Y\text{ such that }(x, y)\in \gamma(t)(V_i)\},\]
which is a finite set. Then we define $P_i$ to be the set of $X$-coordinates  of translations of points in $P_{s, t}$ to the base of $T_i$ 
\[\{x\in X: x=\bar{\alpha}(-m)(z)\text{ for a } z\in P_{s, t}\text{ and  } s, t\in T_i\text{ such that }t=(m, h)\text{ for an }h\},\]
which is still a finite set. Now, for each $i\in I$ define an open set
\[W_i=V_i\setminus \bigcup\{\{x\}\times Y: x\in P_i\}.\]
Note that each $W_i$ is non-empty because 
$\mu(\bigcup\{\{x\}\times Y: x\in P_i\})=0$ for any $\mu\in M_{\Z\times \Gamma}(X\times Y)$. In addition, from this construction, the collection $\CC=\{(S_i, W_i): i\in I\}$  and $\CD=\{(T_i, W_i): i\in I\}$ are open castles in the sense of Definition \ref{defn: tower for partial action}. Then Lemma \ref{lem: extendable-partial action} implies that $\CC$ is $K$-extendable to $\CD$. Moreover, note that for any $i\in I$ and $t\in T_i$, the level $\gamma(t)(W_i)\subseteq \gamma(t)(V_i)\subseteq V$ holds for some $V\in \CV$. Finally, observe that
\[(X\times Y)\setminus \bigsqcup_{i\in I, t\in S_i}\gamma(t)(W_i
)\subseteq R\cup\bigcup_{i\in I}\bigcup_{t\in T_i,x\in P_i}(\{\bar{\alpha}(t)(x)\}\times Y),\]
which implies that 
\[\sup_{\mu\in M_{\Z\times \Gamma}(X\times Y)}\mu((X\times Y)\setminus \bigsqcup_{i\in I, t\in S_i}\gamma(t)(W_i
))<\varepsilon=\inf_{\mu\in M_{\Z\times \Gamma}(X\times Y)}\mu(O).\]
Then Proposition \ref{prop: comparison} shows that $(X\times Y)\setminus \bigsqcup_{i\in I, t\in S_i}\gamma(t)(W_i
)\prec O$ because $\bar\alpha\times\beta$ indeed has the dynamical strict comparison by \cite[Theorem 9.2]{D}.
\end{proof}
    
\begin{thm}
    Let $\beta: \Gamma\curvearrowright Y$ be a minimal free action of an amenable group $G$ on a compact metrizable space $Y$ with $\dim(Y)<\infty$. Then $A=C(Y)\rtimes_r \Gamma$ is $\CZ$-stable if and only if $A$ has a minimal topologically principal second countable almost elementary groupoid model.
\end{thm}
\begin{proof}
    The ``if'' part was proved in \cite{M-W}. For the ``only if'' part, it was shown in \cite[Corollary B]{Na2024} that $\bar{\alpha}\times \beta:\Z\times \Gamma\curvearrowright X\times Y$ is almost finite. Then, Proposition \ref{prop: AE permanance} implies that the transformation groupoid $\CG$ of the partial action $\alpha\times \beta$ is almost elementary. Since $A\simeq A\otimes \CZ\simeq C^*_r(\CG)$, this completes the proof.
\end{proof}

On the other hand, it is still not clear to the authors whether the Jiang-Su algebra $\CZ$ can be written as a crossed product of a dynamical system. It is well known from the Pimsner-Voiculescu exact sequence that $\CZ$ cannot be written as the crossed product of a $\Z$-system. Using more advanced $K$-theoretic tools, we further eliminate most of the usual groups as candidates.

\begin{thm}\label{thm: Jiang Su crossed product}
    Suppose $\CZ\cong C(X)\rtimes_r \Gamma$. Then $\Gamma$ has to be discrete, torsion-free, amenable, and rationally acyclic (that is, all the group homology groups $H_i(B\Gamma; \mathbb{Q})$ with rational coefficients vanish for $i > 0$, where $B\Gamma$ is the classifying space of $\Gamma$).
\end{thm}
\begin{proof}
	The discreteness of $\Gamma$, as well as the compactness of $X$, follows from the unitality of $\CZ$. 
	The group $\Gamma$ must also be torsion-free, because any torsion element $g \in \Gamma$ would give rise to a projection 
	\[
		\frac{1}{|\langle g \rangle|} \sum_{h \in \langle g \rangle} u_h \in C(X)\rtimes_r \Gamma 
	\]
	that is neither $0$ or $1$. 
	
	Since $\CZ$ is nuclear, by \cite[Theorem 5.8]{AD}, the action $\Gamma \curvearrowright X$ is amenable in the sense of \cite[Definition 2.1]{AD}. 
	On the other hand, by restricting the unique trace $\tau$ of $\CZ$ to $C(X)$ (henceforth still denoted by $\tau$), we obtain a $\Gamma$-invariant measure on $X$, whence by \cite[Example 2.7(2)]{AD}, the group $\Gamma$ is amenable. 
	
	It remains to show $\Gamma$ is rationally acyclic. Using the Chern character map 
	\[
		K_i (B\Gamma) \otimes \mathbb{Q}  \xrightarrow{\cong} \bigoplus_{j = i \text{ mod } 2} H_j(B\Gamma; \mathbb{Q}) \; ,
	\]
	it suffices to show the left-hand side agrees with $K_i (\{\mathrm{pt}\}) \otimes \mathbb{Q}$, or equivalently, to show $K_i (B\Gamma) \otimes \mathbb{R} \cong K_i (\{\mathrm{pt}\}) \otimes \mathbb{R}$ as real vector spaces. 
	Now as $\Gamma$ is amenable, by \cite[Corollary 9.2]{HK}, the Baum-Connes assembly map (with the coefficient $\Gamma$-$C^*$-algebra $C(X)$) 
	\[
		KK_*^\Gamma (\underline{E}\Gamma, C(X)) \to K_*(C(X) \rtimes_r \Gamma) 
	\]
	is an isomorphism, where $\underline{E}\Gamma$ is the universal space for proper actions by $\Gamma$. Since $\Gamma$ is torsion-free, $\underline{E}\Gamma$ agrees with ${E}\Gamma$, the universal space for proper and free actions by $\Gamma$, or in other words, the universal cover of $B\Gamma$. 
	Applying the machinery of equivariant $KK$-theory with real coefficients developed in \cite{AAS1} and \cite{AAS2}, we obtain isomorphisms of real vector spaces (where scalar multiplication is given by taking Kasparov products with $KK_{\mathbb{R}}(\mathbb{C}, \mathbb{C}) \cong \mathbb{R}$)
	\[
		KK^\Gamma_{\mathbb{R},*} ({E}\Gamma, C(X)) \cong K_{\mathbb{R}, *}(C(X) \rtimes_r \Gamma) \cong K_{\mathbb{R}, *}(\CZ) \cong K_{\mathbb{R}, *}(\mathbb{C}) \cong  
		\begin{cases}
		\mathbb{R}, & i = 0 \\
		0 , & i = 1 
		\end{cases}
		\; .
	\]
	Now the main reason we apply the machinery of equivariant $KK$-theory with real coefficients is the fact that the invariant trace $\tau$ on $C(X)$ gives rise to an element $[\tau] \in KK^\Gamma_{\mathbb{R}} (C(X), \mathbb{C})$, which is a right-inverse to the element $[\iota] \in KK^\Gamma_{\mathbb{R}} ( \mathbb{C}, C(X))$ given by the unital embedding $\iota \colon \mathbb{C} \to C(X)$, as $\tau \circ \iota = \mathrm{id}_{\mathbb{C}}$. It follows by taking Kasparov products that $KK^\Gamma_{\mathbb{R}, *} ({E}\Gamma, \mathbb{C})$ is a direct summand of $KK^\Gamma_{\mathbb{R}, *} ({E}\Gamma, C(X))$ as real vector spaces, and hence we also have $KK^\Gamma_{\mathbb{R}, *} ({E}\Gamma, \mathbb{C}) \cong K_{\mathbb{R}, *}(\mathbb{C}) \cong K_i (\{\mathrm{pt}\}) \otimes \mathbb{R}$. 
	On the other hand, by \cite[Equation 6.1]{AAS2}, we have $KK^\Gamma_{\mathbb{R}, *} ({E}\Gamma, \mathbb{C}) \cong K_*(B\Gamma) \otimes \mathbb{R}$ as real vector spaces. This completes the proof. 
\end{proof}

\begin{cor}
    Suppose $\CZ\cong C(X)\rtimes_r \Gamma$. Then $\Gamma$ cannot be of type F (that is, $B\Gamma$ cannot be homotopy equivalent to a finite CW complex).
\end{cor}
  
\begin{proof}
    Suppose on the contrary that $B\Gamma$ has a finite CW complex model. Then its Euler characteristic $\chi(\Gamma)$ is well-defined as the signed sum of the Betti numbers, and it agrees with the $\ell^2$-Euler characteristic $\chi_{(2)}(\Gamma)$, that is, the signed sum of the $\ell^2$-Betti numbers (see, e.g., \cite[Theorem~3.19]{Kammeyer2019Introductiona}). By \Cref{thm: Jiang Su crossed product}, $\Gamma$ is infinite amenable, and thus all the $\ell^2$-Betti numbers of $\Gamma$ vanish by a result of Cheeger and Gromov \cite[Theorem~0.2]{CheegerGromov1986L2cohomology}. It follows that $0 = \chi_{(2)}(\Gamma) = \chi(\Gamma)$. Again by \Cref{thm: Jiang Su crossed product}, $\Gamma$ is rationally acyclic, whence all but the $0$-th Betti numbers vanish. Removing these zero terms from the equation $\chi(\Gamma) = \chi_{(2)}(\Gamma)$, we see that the $0$-th Betti number also vanishes, which is absurd, since the retraction of $B\Gamma$ to a singleton implies $H_0(B \Gamma) \not= 0$. 
\end{proof}

\section{Almost elementariness and pure infiniteness for groupoids}
In this section, we first establish the equivalence between almost elementariness and pure infiniteness for groupoids whenever the unit space is compact and there is no groupoid invariant probability Borel measure on the unit space. We first recall the definition of the pure infiniteness defined in \cite{Ma2022Purely} for general groupoids.

\begin{defn}\label{defn: pure inf for groupoid}
    A groupoid $\CG$ is said to be \textit{purely infinite} if $O_1\precsim_{\CG, 2} O_2$ holds in the sense of Definition \ref{defn: groupoid subequivalence} for all non-empty open sets $O_1, O_2$ in $\GU$ satisfying $O_1\subseteq r(\CG O_2)$
\end{defn}

In the minimal case, it follows from  \cite[Theorem 5.1]{Ma2022Purely} that a groupoid $\CG$ is purely infinite if and only if $\CG$ has \textit{groupoid strict comparison} in the sense of \cite[Definition 3.4]{Ma2022Purely} and $M(\CG)=\emptyset$, if and only if $O_1\precsim_\CG O_2$ holds for any non-empty open sets in $\GU$.

As previously mentioned in the introduction, the property of almost elementariness for a groupoid $\CG$ has been established in \cite{M-W} as a crucial property that guarantees the (tracial) $\CZ$-stability of $C^*_r(\CG)$. In the following, we will demonstrate the equivalence between pure infiniteness and almost elementariness for specific groupoids.

\begin{prop}\label{prop: pure inf equal AE}
    Let $\CG$ be a minimal topologically principal groupoid such that $\GU$ is compact. Then $\CG$ is purely infinite if and only if $\CG$ is almost elementary and $M(\CG)=\emptyset$. 
\end{prop}
\begin{proof}
Suppose $\CG$ is almost elementary.  Then \cite[Corollary 8.14]{M-W} implies that $\CG$ has groupoid strict comparison. Now in the case that $M(\CG)=\emptyset$, By \cite[Theorem 5.1]{Ma2022Purely}, one has that $\CG$ is purely infinite.

For the converse, let $K$ be a compact set in $\CG$, $O$ open set in $\GU$, and $\CV$ a finite open cover of $\GU$. Then one can find finitely many open bisections $B_1,\dots, B_n$ such that $K\subseteq \bigcup_{i=1}^nB_i$. Choose a minimal non-empty member in $\{\bigcap_{i\in I}s(B_i): I\subseteq [1, n]\cap \N\}$ with respect to the order ``$\subseteq$'', say $U=\bigcap_{i\in I_0}s(B_i)$. This implies that $U\cap s(B_i)=\emptyset$ for any $i\notin I_0$. Because $\CG$ is topological principal, one chooses a unit $u\in U$ with a trivial stabilizer, which implies that $\{r(B_iu): i\in I_0\}$ contains distinct units. This allows to find open bisections $D_i$ for $i\in I_0$  such that 
\begin{enumerate}[label=(\roman*)]    \item for any $i\in I_0$, the source $s(D_i)=C$ for an open set $C$ in $\GU$ with $u\in C\subseteq \overline{C}\subseteq U$.
    \item $r(D_i)\cap r(D_j)=\emptyset$ for any $i\neq j\in I_0$.
    \item $s(D_i)\cap s(B_j)=\emptyset$ for any $i\in I_0$ and any $j\notin I_0$.
    \item $\{C, r(D_i): i\in I_0\}$ is finer than $\CV$.
\end{enumerate}
Now for $i, j\in I_0$, define a bisection $D_{i, j}=D_iD^{-1}_j$. In addition, we define castles $\CC=\{C\}$ and $\CD=\{D_{i, j}, D_i, D^{-1}_i, C: i, j\in I_0\}$. We  write $B=\bigcup_{i=1}^nB_i$ for simplicity and observe that $\CC$ is $B$-extendable to $\CD$ by our construction of $\CC$ and $\CD$. Finally, because $\CG$ is minimal and purely infinite, one has $\GU\setminus \bigcup\CUU=\GU\setminus C\prec_\CG O$ by \cite[Theorem 5.1]{Ma2022Purely} and we have verified the almost elementariness for $\CG$.
\end{proof}

In \cite{Sp}, Spielberg introduced groupoid models for (non-unital) Kirchberg algebras satisfying the UCT. Then this class of groupoids has been refined in \cite{C-F-H} to cover all unital Kirchberg algebras satisfying the UCT by considering the restriction $\CG_D$ of Spielberg's groupoids $\CG$ to some compact open set $D\subseteq \GU$.

We now show all of Spielberg's groupoids $\CG$ in \cite{Sp} and the restrictions $\CG_D$ in \cite{C-F-H} are minimal topologically principal and purely infinite. In particular, these $\CG_D$'s are almost elementary by Proposition \ref{prop: pure inf equal AE}. This implies that every unital Kirchberg algebras satisfying the UCT has an almost elementary groupoid model.
We refer to \cite[Section 2]{Sp} for the original construction of the groupoid based on a mixture of a 1-graph and two product 2-graphs and we only review the necessary details for our purpose. 

First, the mixture graph described in \cite[Section 2]{Sp}, denoted by $\Omega$, consists of three parts, $E_i\times F_i$ for $i=0, 1$ and $D$, that satisfy the following: 
\begin{enumerate}[label=(\roman*)]    \item  $E_0, E_1, F_0, F_1$ are irreducible directed graphs and $D$ is an explicit irreducible directed graph provided in Figure \ref{figure: D} below. See also \cite[Figure 1]{Sp}.
    \item For $i=0, 1$, denote by $L_i$, respectively, $M_i$ the set of vertices in $E^0_i$, respectively $F^0_i$, emitting infinitely many edges. All $L_i$ and $M_i$ are assumed non-empty. Also, fix distinguished vertices $v_i\in L_i$ and $w_i\in M_i$.
    \item Attach the 2-graphs $E_i\times F_i$ to the graph $D$ by identifying the vertex $u_i\in D^0$ with $(v_i, w_i)\in E^0_i\times F^0_i$ for each $i=1,2$.
\end{enumerate}

\begin{figure}[h]
    \centering
    \tikzset{
        midarrow/.style={
            postaction={
                decorate,
                decoration={
                    markings,
                    mark=at position 0.5 with {\arrow{>}}
                }
            }
        }
    }
    \begin{tikzpicture}[>=latex, scale=1.2, thick]
        \node (a0) at (0, 2.5) {\(a_0\)};
        \node (u0) at (-3.5, 0) {\(u_0\)};
        \node (u1) at (3.5, 0) {\(u_1\)};
        \node (a1) at (0, -2.5) {\(a_1\)};

        \draw[midarrow] (u0) -- (a0) node[midway, above left, inner sep=1pt] {\(\alpha_0\)};
        \draw[midarrow] (a0) -- (u1) node[midway, above right, inner sep=1pt] {\(\delta_0\)};
        \draw[midarrow] (a1) -- (u0) node[midway, below left, inner sep=1pt] {\(\delta_1\)};
        \draw[midarrow] (u1) -- (a1) node[midway, below right, inner sep=1pt] {\(\alpha_1\)};

        \draw[midarrow] (u0) to[bend left=35] node[midway, above, inner sep=2pt] {\(\varepsilon_0\)} (u1);
        \draw[midarrow] (u1) to[bend left=35] node[midway, below, inner sep=2pt] {\(\varepsilon_1\)} (u0);

        \draw[->] (a0) .. controls (-2.5, 3.8) and (-2.5, 2.5) .. (a0);
        \node at (-1.2, 3.3) {\(\beta_0\)};

        \draw[->] (a0) .. controls (2.5, 3.8) and (2.5, 2.5) .. (a0);
        \node at (1.2, 3.3) {\(\gamma_0\)};

        \draw[->] (a1) .. controls (-2.5, -3.8) and (-2.5, -2.5) .. (a1);
        \node at (-1.2, -3.3) {\(\gamma_1\)};

        \draw[->] (a1) .. controls (2.5, -3.8) and (2.5, -2.5) .. (a1);
        \node at (1.2, -3.3) {\(\beta_1\)};

    \end{tikzpicture}
    
    \caption{The graph $D$}\label{figure: D}
\end{figure}

\begin{rmk}
    We remark that the irreducibility of a directed graph $G$ means that for any vertices $u,v\in G$, there is a directed path connecting $u$ and $v$. Some authors use the term \textit{strong connectedness} instead of irreducibility. 
\end{rmk}

A \textit{vertex} of $\Omega$ refers to an element of $\bigcup_{i=0, 1}(E^0_i\times F^0_i)\cup D^0$ with the identification of $u_i$ with $(v_i, w_i)$. An \textit{edge} of $\Omega$ means an element in 
\[(\bigcup_{i=0, 1}(E^1_i\times F^0_i)\cup (E^0_i\times F^1_i))\cup D^1.\]

The groupoid $\CG$ constructed from $\Omega$ is similar to the usual graph groupoid. In particular, the unit space $\GU$ consists of specific
paths in $\Omega$. We recall the following definition in \cite{Sp}.

\begin{defn}\cite[Definition 2.3]{Sp}
    A \textit{finite path element of type $D$} is a finite directed path in $D$ with non-zero length. An \textit{infinite path element of type $D$} is either an infinite directed path in $D$ or a finite path element of type $D$ that ends at $u_0$ or $u_1$. A \textit{finite path element of type $(E_i, F_i)$} is an ordered pair $(p, q)\in E^*\times F^*$ such that at least one of $p, q$ is of positive length. An \textit{infinite path element of type $(E_i, F_i)$} is an ordered pair $(p, q)$, where $p$, respectively $q$, is either an infinite path  or a finite path terminating in $L_i$, respectively in $M_i$, in $E_i$, respectively $F_i$, and $(p, q)$ are not both of length zero.

    Also define the \textit{origin} and the \textit{terminus} for a path element $(p, q)$ of type $(E_i, F_i)$ by, respectively, $o(p, q)=(o(p), o(q))$ and $t(p, q)=(t(p), t(q))$. We also say $(p, q)$ \textit{extends} $(p', q')$ if $p$ extends $p'$ and $q$ extends $q'$ in the usual sense.
\end{defn}

Using path elements of various types, one can build paths in the following way.

\begin{defn}\cite{Sp}\label{defn: paths}
    A \textit{finite path} is either a vertex in $\Omega$, or a finite string $\mu_1\cdots\mu_k$ of finite path element such that  such that 
    \begin{enumerate}[label=(\roman*)]        \item $t(\mu_i)=o(\mu_{i+1})$, and
    \item $\mu_i$ and $\mu_{i+1}$ are of different types.
    \end{enumerate}
An \textit{infinite path} is either a vertex in $L_i\times M_i$, an infinite string of finite path elements satisfying (i) and (ii) above, or a finite sequence $\mu_1\dots\mu_{k+1}$ satisfying
\begin{enumerate}
    \item[(iii)] $\mu_1\dots\mu_k$ is a finite path in the sense above.
    \item[(iv)] $\mu_{k+1}$ is an infinite path element.
    \item[(v)] (i) and (ii) above hold.
\end{enumerate}
\end{defn}

Denote by $X$ the set of all infinite paths as the unit space for the groupoid that will be defined. Also write $\mu\preceq \nu$ for the situation that $\nu$ extends $\mu$. Then we recall the definition of topology on $X$.

\begin{defn}\cite[Definition 2.8]{Sp}
    Let $\mu=\mu_1\dots\mu_k$ be a finite path in the sense of Definition \ref{defn: paths}. define
    \[Z(\mu)=\{\sigma=\sigma_1\sigma_2\dots\in X: \sigma_i=\mu_i\ \text{for }i<k,\text{ and } \mu_k\preceq \sigma_k\}\]
    and define $V(\mu)\subseteq Z(\mu)$ by either $V(\mu)=Z(\mu)$ when $t(\mu)\notin \{u_0, u_1\}$, or $V(\mu)=(Z(\mu)\setminus Z(\mu\alpha_i))\setminus Z(\mu\varepsilon_i)$ if $t(\mu)=u_i$, in which $\alpha_i$ and $\varepsilon_i$ are particular edges in $D$ (see \cite[Figure 1]{Sp}). Moreover,  in the case that $t(\mu)=(y, z)\in E^0_i\times F^0_i$, for any finite sets $B\subseteq E^1(y)$ and $C\subseteq F^1(z)$ such that if $y\notin L_i$ then $B=\emptyset$ and if $z\notin M_i$ then $C=\emptyset$, one defines
    \[V(\mu; B, C)=(V(\mu)\setminus \bigcup_{e\in B}Z(\mu(e, z)))\setminus \bigcup_{f\in C}Z(\mu(y,f)),\]
    in which $\mu(e, z)$ and $\mu(y, f)$ are concatenations of $\mu$ with edges $(e, z)$ and $(y, f)$, respectively.
\end{defn}

These sets are expected to generate the topology on $X$ and to induce bisections on the desired groupoid. Indeed, It was shown in \cite[Lemma 2.14]{Sp} that the collection $\CB$ of all possible $Z(\mu)$ and $V(\mu; B, C)$ forms a base for a locally compact metrizable topology on $X$. 

\begin{defn}
    Define the \textit{length function} $\ell :\{\text{finite path elements}\}\to \N^2$ by 
    \[\ell(\mu)=\begin{cases}
    (0, 0), & \mu\text{ is a vertex,}\\
(\ell(p), \ell(p)), & \mu\text{ is a finite path element of type }D,\\
(\ell(p), \ell(q)), & \mu\text{ is a finite path element of type }(E_i, F_i).
\end{cases}
\]
Then extend the definition of $\ell$ to finite path $\mu=\mu_1\dots\mu_k$ by \[\ell(\mu)=\sum_{i=1}^k\ell(\mu_k).\]
\end{defn}
Let $\CG$ be the set of triples $(x, n, y)\in X\times \Z^2\times X$ such that there exists $z\in X$ and decompositions $x=\mu z$ and $y=\nu z$ with $\ell(\mu)-\ell(\nu)=n$. 

\begin{defn}\cite{Sp}\label{defn: groupoid for Kirchberg algebra}
    $\CG$ defined above is a locally compact ample Hausdorff \'{e}tale groupoid when it is equipped with the multiplication $(x, n, y)(y, m, z)=(x, n+m, z)$ and the inverse operation $(x, n, y)^{-1}=(y, -n, x)$. The topology on $\CG$ is generated by the following compact open bisections (as a base):
    \begin{align*}
 U(\mu_1, \mu_2)&=(Z(\mu_1)\times \{\ell(\mu_1)-\ell(\mu_2)\}\times Z(\mu_2))\cap \CG\\
 U_0(\mu_1, \mu_2; B, C)&=(V(\mu_1; B, C)\times \{\ell(\mu_1)-\ell(\mu_2)\}\times V(\mu_2; B, C))\cap \CG    \end{align*}
 where $\mu_1$ and $\mu_2$ are finite paths with $t(\mu_1)=t(\mu_2)$.
\end{defn}

The final result we will recall from \cite{Sp} is the following regarding the basic fundamental properties of $\CG$.

\begin{lem}\cite[Lemma 2.18]{Sp}\label{lem: basic property for groupoid for Kirchberg algebras}
    $\CG$ is minimal, topological principal, amenable, and locally contractive.
\end{lem}

We remark that it was proved in \cite[Theorem 5.1, Corollary 5.7]{Ma2022Purely} that pure infiniteness implies the local contractivity for any minimal ample groupoid. The converse is not true in general (see \cite[Theorem D]{Ma2022Purely} ) but is not known for the minimal ample case. However, in the rest of this section, we establish the pure infiniteness for $\CG$ in Definition \ref{defn: groupoid for Kirchberg algebra} directly and thus strengthen Lemma \ref{lem: basic property for groupoid for Kirchberg algebras}. We begin with the following lemma.

\begin{prop}\label{prop: pure inf non-unital}
    Let $\CG$ be the groupoid defined in Definition \ref{defn: groupoid for Kirchberg algebra}. Then $\CG$ is purely infinite.
\end{prop}
\begin{proof}
    Using \cite[Theorem 5.1]{Ma2022Purely}, it suffices to show $K\prec_\CG O$ holds for any compact set $K$ and non-empty open set $O$ in $X=\GU$. Note that for any finite path $\mu$ with $t(\mu)=(y, z)\in E^0_i\times F^0_i$, one has $V(\mu; B, C)\subseteq Z(\mu)$ for any possible $B\subseteq E^1(y)$ and $C\subseteq F^1(z)$. Moreover, extend a finite path if necessary, one may assume there is a finite path $\mu$ with $t(\mu)\in E^0_i\times F^0_i$ for some $i=0, 1$ such that $V(\mu; B, C)\subseteq O$. Now we show $K\prec_\CG V(\mu; B, C)$.

    Choose finitely many open sets $Z(\sigma_k)$ for $k=1,\dots, m$ such that 
    $K\subseteq  \bigcup_{k=1}^mZ(\sigma_k).$ Then choose $m$ disjoint open sets of form $V(\gamma_l; B'_l, C'_l)$ contained in $V(\mu; B, C)$ for $l=1,\dots, m+n$. Note that $\Omega$ is irreducible as a directed graph because $E_i, F_i$ ($i=0,1$) and $D$ are.
Therefore, for $l=1, \dots, m$, one extends $\gamma_l$ to $\gamma'_l$ in a way such that $t(\gamma'_l)=t(\sigma_l)$ and $Z(\gamma'_l)\subseteq V(\gamma_l; B'_l, C'_l)$. Now using bisections $U(\gamma'_l, \sigma_l)$ for $l=1,\dots, m$  in Definition \ref{defn: groupoid for Kirchberg algebra}, one has 
\[K\subseteq\bigcup_{l=1}^ms(U(\gamma'_l, \sigma_k))\]
and
\[\bigsqcup_{l=1}^mr(U(\gamma'_k, \sigma_k))\subseteq \bigsqcup_{l=1}^{m}V(\gamma_l; B'_l, C'_l)\subseteq V(\mu; B, C),\]
which implies that $K\prec_\CG V(\mu; B, C)$. This shows that $\CG$ is purely infinite.
\end{proof}

Let $\CG$ be the groupoid in Definition \ref{defn: groupoid for Kirchberg algebra}. Fix an vertex $(y,z)\in E^0_0\times F^0_0$ and denote by $\mu=(y, z)$. It was proved in \cite{C-F-H} that $K_*(C^*_r(\CG))=K_*(C^*_r(\CG|_{Z(\mu)}))$, where $\CG|_{Z(\mu)}$ is the restriction of $\CG$ on the compact open set $Z(\mu)$. It is not hard to see directly that $\CG|_{Z(\mu)}$ is also minimal, topologically principal, and amenable. We now show $\CG|_{Z(\mu)}$ is purely infinite as well. Note that for any $S\subseteq \CG$, the restriction $S|_{Z(\mu)}$ of $S$ on $Z(\mu)$ is exactly the set $USU$ in $\CG$. 

\begin{prop}\label{prop: pure inf unital}
The groupoid $\CH=\CG|_{Z(\mu)}$ is purely infinite and thus almost elementary.
\end{prop}
\begin{proof}
    Note that $\HU=Z(\mu)$, which is a compact set. To show $\CH$ is purely infinite, it suffices to show $\HU\prec_\CH O$ for any non-empty open set $O$ in $\HU$.
    
 Proposition \ref{prop: pure inf non-unital} implies that $\HU\prec_\CG O$ in $\CG$ and thus there are open bisections $A_1,\dots, A_n$ in $\CG$ such that $\HU\subseteq \bigcup_{i=1}^ns(A_i)$ and $\bigsqcup_{i=1}^nr(A_i)\subseteq O$. 
 Now, define open bisections $B_i=UA_iU$ in $\CH$ and observe that $\HU=U\HU U$ and $O=UOU$. Therefore, one has $\HU\subseteq \bigcup_{i=1}^ns(B_i)$ and $\bigsqcup_{i=1}^nr(B_i)\subseteq O$, which means $\HU\prec_\CH O$.

 Thus, $\CH$ is purely infinite and thus almost elementary by Proposition \ref{prop: pure inf equal AE}.
  \end{proof}

 Using Proposition \ref{prop: pure inf non-unital}, \ref{prop: pure inf unital} and \cite{Sp} and \cite{C-F-H}, we have arrived at the following conclusion.

\begin{thm}\label{thm: groupoid model for kirchberg}
    Every Kirchberg algebra satisfying the UCT has a minimal topologically principal, purely infinite groupoid model. In particular, every unital Kirchberg algebra satisfying the UCT has a minimal, topologically principal, topologically amenable, almost elementary groupoid model.
\end{thm}

\begin{rmk}
   Most recently, it has been shown in \cite{BK2026} (see also \cite{Wu2025Cstar}) any unital Kichberg algebra has a minimal topologically principal topologically amenable groupoid model $\CG$ arising from partial dynamics. Moreover, as shown in the proof of \cite[Theorem 7.6]{BK2026}, the $\CG$ can be even chosen to have \textit{paradoxical comparison} in the sense of \cite[Definition 3.6]{Ma2022Purely}, which implies that $\CG$ is purely infinite by \cite[Theorem 5.1]{Ma2022Purely}. Therefore, Proposition \ref{prop: pure inf equal AE} shows that $\CG$ is almost elementary. Thus, this process provides a new type of minimal, topologically principal, topologically amenable, almost elementary groupoid model for all unital Kirchberg algebras satisfying the UCT.
\end{rmk}

\section*{Ackowledgement}
This research is supported by National Key R\&D Program of China 2022YFA100700. 
The first author is supported by NSFC No.~12571133. 
The second author is partially supported by NSFC Key Program No.~12231005. 

The authors would like to thank N. Christopher Phillips, Aaron Tikuisis, and Xiaolei Wu for very helpful discussions.

 \bibliographystyle{alpha}
 \bibliography{AE_model_for_ssa_algebras}

@article{EGLN,
    author = {Elliott, George A. and Gong, Guihua and Lin, Huaxin and Niu, Zhuang},
    title = {On the classification of simple amenable {C}*-algebras with finite decomposition rank, {II}},
    journal = {J. Noncommut. Geom.},
    volume = {19},
    number = {1},
    pages = {73--104},
    year = {2025},
    mrnumber = {MR4860189}
}

@Article{D,
  author   = {Kerr, David},
  title    = {Dimension, comparison, and almost finiteness},
  journal  = {J. Eur. Math. Soc. (JEMS)},
  year     = {2020},
  volume   = {22},
  number   = {11},
  pages    = {3697--3745},
  issn     = {1435-9855},
  doi      = {10.4171/jems/995},
  fjournal = {Journal of the European Mathematical Society (JEMS)},
  mrclass  = {Prelim},
  mrnumber = {4167017},
  url      = {https://doi.org/10.4171/jems/995},
}

@Article{K-S,
  author   = {Kerr, David and Szab\'{o}, G\'{a}bor},
  title    = {Almost finiteness and the small boundary property},
  journal  = {Comm. Math. Phys.},
  year     = {2020},
  volume   = {374},
  number   = {1},
  pages    = {1--31},
  issn     = {0010-3616},
  doi      = {10.1007/s00220-019-03519-z},
  fjournal = {Communications in Mathematical Physics},
  mrclass  = {22D25 (19K99 37A25 37A55)},
  mrnumber = {4066584},
  sortkey  = {KerrSzabo2020Almost},
  url      = {https://doi.org/10.1007/s00220-019-03519-z},
}

@Article{Li,
  author     = {Li, Xin},
  title      = {Every classifiable simple $C^*$-algebra has a Cartan subalgebra},
  journal    = {Invent. Math.},
  year       = {2020},
  volume     = {219},
  number     = {2},
  pages      = {653--699},
  issn       = {0020-9910},
  doi        = {10.1007/s00222-019-00914-0},
  fjournal   = {Inventiones Mathematicae},
  mrclass    = {46L35 (22A22)},
  mrnumber   = {4054809},
  mrreviewer = {James Gabe},
  url        = {https://doi.org/10.1007/s00222-019-00914-0},
}

@article{Ma2022Purely,
      author       = {Xin Ma},
      title        = {Purely Infinite Locally Compact Hausdorff étale Groupoids and Their C*-algebras},
      journal      = {Int. Math. Res. Not.(IMRN)},
      volume       = {2022},
      number       = {11},
      pages        = {8420--8471},
      year         = {2022},
      doi          = {10.1093/imrn/rnaa360},
      url          = {https://doi.org/10.1093/imrn/rnaa360},
      published    = {30 January 2021},
      received     = {25 May 2020},
}

@Article{Matui2,
  author     = {Matui, Hiroki},
  title      = {Topological full groups of one-sided shifts of finite type},
  journal    = {J. Reine Angew. Math.},
  year       = {2015},
  volume     = {705},
  pages      = {35--84},
  issn       = {0075-4102},
  doi        = {10.1515/crelle-2013-0041},
  fjournal   = {Journal f\"{u}r die Reine und Angewandte Mathematik. [Crelle's Journal]},
  mrclass    = {22A22 (37B10)},
  mrnumber   = {3377390},
  mrreviewer = {Thomas Ward},
  url        = {https://doi.org/10.1515/crelle-2013-0041},
}

@Article{Matui,
  author     = {Matui, Hiroki},
  title      = {Homology and topological full groups of \'{e}tale groupoids on totally disconnected spaces},
  journal    = {Proc. Lond. Math. Soc. (3)},
  year       = {2012},
  volume     = {104},
  number     = {1},
  pages      = {27--56},
  issn       = {0024-6115},
  doi        = {10.1112/plms/pdr029},
  fjournal   = {Proceedings of the London Mathematical Society. Third Series},
  mrclass    = {19D55 (22A22 46L55)},
  mrnumber   = {2876963},
  mrreviewer = {Hiroyuki Osaka},
  url        = {https://doi.org/10.1112/plms/pdr029},
}

@article{Na2024,
    author = {Naryshkin, Petr},
    title = {Group extensions preserve almost finiteness},
    journal = {J. Funct. Anal.},
    volume = {286},
    number = {7},
    pages = {110348},
    year = {2024},
    doi = {10.1016/j.jfa.2024.110348},
    zbl = {1544.46055}
}

@Article{Ne,
  author     = {Nekrashevych, Volodymyr},
  title      = {Simple groups of dynamical origin},
  journal    = {Ergodic Theory Dynam. Systems},
  year       = {2019},
  volume     = {39},
  number     = {3},
  pages      = {707--732},
  issn       = {0143-3857},
  doi        = {10.1017/etds.2017.47},
  fjournal   = {Ergodic Theory and Dynamical Systems},
  mrclass    = {22D40 (22A22 37B05)},
  mrnumber   = {3904185},
  mrreviewer = {Ahmad Zainy Al-Yasry},
  url        = {https://doi.org/10.1017/etds.2017.47},
}

@Book{Renault,
  title      = {A groupoid approach to {$C^{\ast} $}-algebras},
  publisher  = {Springer, Berlin},
  year       = {1980},
  author     = {Renault, Jean},
  volume     = {793},
  series     = {Lecture Notes in Mathematics},
  isbn       = {3-540-09977-8},
  mrclass    = {46Lxx (22D25 22D40)},
  mrnumber   = {584266},
  mrreviewer = {A. K. Seda},
  pages      = {ii+160},
}

@Unpublished{Sims,
  author = {Sims, Aidan},
  title  = {{H}ausdorff \'{e}tale groupoids and their ${C}^{*}$-algebras},
  note   = {arXiv:1710.10897},
  year   = {2017},
}

@Article{TWW,
  author     = {Tikuisis, Aaron and White, Stuart and Winter, Wilhelm},
  title      = {Quasidiagonality of nuclear {$C^\ast$}-algebras},
  journal    = {Ann. of Math. (2)},
  year       = {2017},
  volume     = {185},
  number     = {1},
  pages      = {229--284},
  issn       = {0003-486X},
  doi        = {10.4007/annals.2017.185.1.4},
  fjournal   = {Annals of Mathematics. Second Series},
  mrclass    = {46L05 (47L40)},
  mrnumber   = {3583354},
  mrreviewer = {Dinesh Jayantilal Karia},
  url        = {https://doi.org/10.4007/annals.2017.185.1.4},
}

@Article{AD,
  author     = {Anantharaman-Delaroche, Claire},
  title      = {Amenability and exactness for dynamical systems and their {$C^\ast$}-algebras},
  journal    = {Trans. Amer. Math. Soc.},
  year       = {2002},
  volume     = {354},
  number     = {10},
  pages      = {4153--4178 (electronic)},
  issn       = {0002-9947},
  coden      = {TAMTAM},
  doi        = {10.1090/S0002-9947-02-02978-1},
  fjournal   = {Transactions of the American Mathematical Society},
  mrclass    = {46L55 (22D25 43A07)},
  mrnumber   = {1926869 (2004e:46082)},
  mrreviewer = {Berndt Brenken},
  url        = {http://dx.doi.org/10.1090/S0002-9947-02-02978-1},
}

@Article{AAS1,
  author        = {Antonini, Paolo and Azzali, Sara and Skandalis, Georges},
  title         = {Bivariant {$K$}-theory with {$\Bbb{R}/\Bbb{Z}$}-coefficients and rho classes of unitary representations},
  journal       = {J. Funct. Anal.},
  year          = {2016},
  volume        = {270},
  number        = {1},
  pages         = {447--481},
  issn          = {0022-1236,1096-0783},
  __markedentry = {[valsvp:6]},
  doi           = {10.1016/j.jfa.2015.06.017},
  fjournal      = {Journal of Functional Analysis},
  mrclass       = {46L80 (19K35)},
  mrnumber      = {3419768},
  mrreviewer    = {Michael\ Frank},
  url           = {https://doi.org/10.1016/j.jfa.2015.06.017},
}

@Article{AAS2,
  author     = {Antonini, Paolo and Azzali, Sara and Skandalis, Georges},
  title      = {The {B}aum-{C}onnes conjecture localised at the unit element of a discrete group},
  journal    = {Compos. Math.},
  year       = {2020},
  volume     = {156},
  number     = {12},
  pages      = {2536--2559},
  issn       = {0010-437X},
  doi        = {10.1112/s0010437x20007502},
  fjournal   = {Compositio Mathematica},
  mrclass    = {19K35 (46L80 58J22)},
  mrnumber   = {4202454},
  mrreviewer = {Sayan Chakraborty},
  url        = {https://doi.org/10.1112/s0010437x20007502},
}

@Article{CETWW,
  author     = {Castillejos, Jorge and Evington, Samuel and Tikuisis, Aaron and White, Stuart and Winter, Wilhelm},
  title      = {Nuclear dimension of simple {$\rm C^*$}-algebras},
  journal    = {Invent. Math.},
  year       = {2021},
  volume     = {224},
  number     = {1},
  pages      = {245--290},
  issn       = {0020-9910,1432-1297},
  doi        = {10.1007/s00222-020-01013-1},
  fjournal   = {Inventiones Mathematicae},
  mrclass    = {46L35 (46L05)},
  mrnumber   = {4228503},
  mrreviewer = {Changguo\ Wei},
  url        = {https://doi.org/10.1007/s00222-020-01013-1},
}

@Article{DPS2,
  author        = {Deeley, Robin J. and Putnam, Ian F. and Strung, Karen R.},
  title         = {Classifiable {$\rm C^*$}-algebras from minimal {$\Bbb Z$}-actions and their orbit-breaking subalgebras},
  journal       = {Math. Ann.},
  year          = {2024},
  volume        = {388},
  number        = {1},
  pages         = {703--729},
  issn          = {0025-5831,1432-1807},
  __markedentry = {[valsvp:6]},
  doi           = {10.1007/s00208-022-02526-1},
  fjournal      = {Mathematische Annalen},
  mrclass       = {37B05 (19K99 37A55 46L35 46L85)},
  mrnumber      = {4693944},
  url           = {https://doi.org/10.1007/s00208-022-02526-1},
}

@Article{DPS1,
  author        = {Deeley, Robin J. and Putnam, Ian F. and Strung, Karen R.},
  title         = {Constructing minimal homeomorphisms on point-like spaces and a dynamical presentation of the {J}iang-{S}u algebra},
  journal       = {J. Reine Angew. Math.},
  year          = {2018},
  volume        = {742},
  pages         = {241--261},
  issn          = {0075-4102,1435-5345},
  __markedentry = {[valsvp:6]},
  doi           = {10.1515/crelle-2015-0091},
  fjournal      = {Journal f\"ur die Reine und Angewandte Mathematik. [Crelle's Journal]},
  mrclass       = {46L40 (37A55 37B05)},
  mrnumber      = {3849627},
  mrreviewer    = {Marius\ V.\ Ionescu},
  url           = {https://doi.org/10.1515/crelle-2015-0091},
}

@Book{GKPT,
  title      = {Crossed products of {$C^*$}-algebras, topological dynamics, and classification},
  publisher  = {Birkh\"{a}user/Springer, Cham},
  year       = {2018},
  author     = {Giordano, Thierry and Kerr, David and Phillips, N. Christopher and Toms, Andrew},
  series     = {Advanced Courses in Mathematics. CRM Barcelona},
  isbn       = {978-3-319-70868-3; 978-3-319-70869-0},
  note       = {Lecture notes based on the course held at the Centre de Recerca Matem\`atica (CRM) Barcelona, June 14--23, 2011, Edited by Francesc Perera},
  doi        = {10.1007/978-3-319-70869-0},
  mrclass    = {46-02 (46L05 46L10 46L40 46L55)},
  mrnumber   = {3837150},
  mrreviewer = {Efton Park},
  pages      = {x+498},
  url        = {https://doi-org.srv-proxy1.library.tamu.edu/10.1007/978-3-319-70869-0},
}

@Article{HK,
  author     = {Higson, Nigel and Kasparov, Gennadi},
  title      = {{$E$}-theory and {$KK$}-theory for groups which act properly and isometrically on {H}ilbert space},
  journal    = {Invent. Math.},
  year       = {2001},
  volume     = {144},
  number     = {1},
  pages      = {23--74},
  issn       = {0020-9910},
  coden      = {INVMBH},
  doi        = {10.1007/s002220000118},
  fjournal   = {Inventiones Mathematicae},
  mrclass    = {19K35 (19L47 46L80)},
  mrnumber   = {1821144 (2002k:19005)},
  mrreviewer = {Emmanuel C. Germain},
  url        = {http://dx.doi.org.proxy.library.vanderbilt.edu/10.1007/s002220000118},
}

@Unpublished{H-W,
  author = {Hirshberg, Ilan and Wu, Jianchao},
  title  = {Long thin covers and nuclear dimension},
  note   = {arXiv:2308.12504},
  year   = {2023},
  url    = {http://arxiv.org/abs/2308.12504},
}

@article{Jos,
    author = {Joseph, Matthieu},
    title = {Amenable wreath products with non almost finite actions of mean dimension zero},
    journal = {Trans. Amer. Math. Soc.},
    volume = {377},
    number = {2},
    pages = {1321--1333},
    year = {2024},
    doi = {10.1090/tran/9061},
    zbl = {1535.37003}
}

@Article{Sp,
  author     = {Spielberg, Jack},
  title      = {Graph-based models for {K}irchberg algebras},
  journal    = {J. Operator Theory},
  year       = {2007},
  volume     = {57},
  number     = {2},
  pages      = {347--374},
  issn       = {0379-4024},
  fjournal   = {Journal of Operator Theory},
  mrclass    = {46L05 (46L35 46L80)},
  mrnumber   = {2329002},
  mrreviewer = {Aidan D. Sims},
}

@Article{C-F-H,
  author        = {Clark, Lisa Orloff and Fletcher, James and an Huef, Astrid},
  title         = {All classifiable {K}irchberg algebras are {$C^*$}-algebras of ample groupoids},
  journal       = {Expo. Math.},
  year          = {2020},
  volume        = {38},
  number        = {4},
  pages         = {559--565},
  issn          = {0723-0869,1878-0792},
  __markedentry = {[valsvp:6]},
  doi           = {10.1016/j.exmath.2019.06.001},
  fjournal      = {Expositiones Mathematicae},
  mrclass       = {46L05 (46L35)},
  mrnumber      = {4177957},
  url           = {https://doi.org/10.1016/j.exmath.2019.06.001},
}

@Unpublished{CGSTW,
  Title                    = {Classifying $*$-homomorphisms {$\rm I$}: Unital simple nulcear {$\rm C^*$}-algebras},
  Author                   = {Carri\'{o}n, Jos\'{e} R. and Gabe, James and Schafhauser, Christopher and Tikuisis, Aaron and White, Stuart},
  Note                     = {arXiv: 2307.06480},
  Year                     = {2023}
}

@Unpublished{M-W,
  Title                    = {Almost elementary \'{e}tale groupoids},
  Author                   = {Ma, Xin and Wu, Jianchao},
  Note                     = {preprint, arXiv:2011.01182},
  Year                     = {2020}
}

@Unpublished{MWfiberwise,
  Title                    = {Fiberwise amenability of \'{e}tale groupoids},
  Author                   = {Ma, Xin and Wu, Jianchao},
  Note                     = {preprint},
  Year                     = {2026}
}

@Book{Exel,
  title         = {Partial dynamical systems, {F}ell bundles and applications},
  publisher     = {American Mathematical Society, Providence, RI},
  year          = {2017},
  author        = {Exel, Ruy},
  volume        = {224},
  series        = {Mathematical Surveys and Monographs},
  isbn          = {978-1-4704-3785-5},
  __markedentry = {[valsvp:]},
  doi           = {10.1090/surv/224},
  mrclass       = {46L55 (16S35 16S40 37A55 46-02 46L45)},
  mrnumber      = {3699795},
  mrreviewer    = {Fernando\ Abadie},
  pages         = {vi+321},
  url           = {https://doi.org/10.1090/surv/224},
}

@Unpublished{GGGKN,
  author = {Gardella, Eusebio and Geffen, Shirly and Gesing, Rafaela and Kopsacheilis, Grigoris and Naryshkin, Petr},
  title  = {Essential Freeness, Allostery and $\mathcal{Z}$-Stability of Crossed Products},
  note   = {arXiv:2405.04343},
  year   = {2024},
  url    = {http://arxiv.org/abs/2308.12504},
}

@Unpublished{BK2026,
  author = {Buss, Alcides and Kranz, Julian},
  title  = {On groupoids beyond partial actions, inner amenability, and models for {Kirchberg} algebras},
  note   = {arXiv:2604.17921},
  year   = {2026},
  url    = {http://arxiv.org/abs/2604.17921},
}

@Article{JiangSu,
  author        = {Jiang, Xinhui and Su, Hongbing},
  title         = {On a simple unital projectionless {$C^*$}-algebra},
  journal       = {Amer. J. Math.},
  year          = {1999},
  volume        = {121},
  number        = {2},
  pages         = {359--413},
  issn          = {0002-9327,1080-6377},
  __markedentry = {[valsvp:6]},
  fjournal      = {American Journal of Mathematics},
  mrclass       = {46L35 (19K35 46L80)},
  mrnumber      = {1680321},
  mrreviewer    = {Vicumpriya\ S.\ Perera},
  url           = {http://muse.jhu.edu/journals/american_journal_of_mathematics/v121/121.2jiang.pdf},
}

@Article{BBWW,
  author        = {Brown, Nathanial and Browne, Sarah L. and Willett, Rufus and Wu, Jianchao},
  title         = {The {UCT} problem for nuclear {$C^*$}-algebras},
  journal       = {Rocky Mountain J. Math.},
  year          = {2022},
  volume        = {52},
  number        = {3},
  pages         = {817--827},
  issn          = {0035-7596,1945-3795},
  __markedentry = {[valsvp:6]},
  doi           = {10.1216/rmj.2022.52.817},
  fjournal      = {The Rocky Mountain Journal of Mathematics},
  mrclass       = {46L35 (19K35)},
  mrnumber      = {4441099},
  mrreviewer    = {Prahlad\ Vaidyanathan},
  url           = {https://doi.org/10.1216/rmj.2022.52.817},
}

@Article{Tu,
  author        = {Tu, Jean-Louis},
  title         = {La conjecture de {B}aum-{C}onnes pour les feuilletages moyennables},
  journal       = {$K$-Theory},
  year          = {1999},
  volume        = {17},
  number        = {3},
  pages         = {215--264},
  issn          = {0920-3036,1573-0514},
  __markedentry = {[valsvp:6]},
  doi           = {10.1023/A:1007744304422},
  fjournal      = {$K$-Theory. An Interdisciplinary Journal for the Development, Application, and Influence of $K$-Theory in the Mathematical Sciences},
  mrclass       = {19K35 (46L85 58J22)},
  mrnumber      = {1703305},
  mrreviewer    = {Hiroshi\ Takai},
  url           = {https://doi.org/10.1023/A:1007744304422},
}

@article{Wu2025CStar,
    author = {Wu, Victor},
    title = {{$C^*$}-algebras associated to directed graphs of groups, and models of {Kirchberg} algebras},
    journal = {J. Funct. Anal.},
    volume = {288},
    number = {3},
    pages = {110740},
    year = {2025},
    doi = {10.1016/j.jfa.2024.110740},
    issn = {0022-1236}
}

@Article{BarlakLi,
  author     = {Barlak, Sel\c{c}uk and Li, Xin},
  title      = {Cartan subalgebras and the {UCT} problem},
  journal    = {Adv. Math.},
  year       = {2017},
  volume     = {316},
  pages      = {748--769},
  issn       = {0001-8708},
  doi        = {10.1016/j.aim.2017.06.024},
  fjournal   = {Advances in Mathematics},
  mrclass    = {46L10},
  mrnumber   = {3672919},
  mrreviewer = {Xiao Chen},
  url        = {https://doi-org.srv-proxy1.library.tamu.edu/10.1016/j.aim.2017.06.024},
}

@Article{BCFS,
  author        = {Brown, Jonathan and Clark, Lisa Orloff and Farthing, Cynthia and Sims, Aidan},
  title         = {Simplicity of algebras associated to \'etale groupoids},
  journal       = {Semigroup Forum},
  year          = {2014},
  volume        = {88},
  number        = {2},
  pages         = {433--452},
  issn          = {0037-1912,1432-2137},
  __markedentry = {[valsvp:6]},
  doi           = {10.1007/s00233-013-9546-z},
  fjournal      = {Semigroup Forum},
  mrclass       = {46L05 (22A22)},
  mrnumber      = {3189105},
  mrreviewer    = {Jean\ N.\ Renault},
  url           = {https://doi.org/10.1007/s00233-013-9546-z},
}

@article{GLN2020a,
    author = {Gong, Guihua and Lin, Huaxin and Niu, Zhuang},
    title = {A classification of finite simple amenable {Z}-stable {$C^*$}-algebras, {I}: {$C^*$}-algebras with generalized tracial rank one},
    journal = {C. R. Math. Rep. Acad. Sci. Canada},
    volume = {42},
    number = {3},
    pages = {63--450},
    year = {2020}
}

@article{GLN2020b,
    author = {Gong, Guihua and Lin, Huaxin and Niu, Zhuang},
    title = {A classification of finite simple amenable {Z}-stable {$C^*$}-algebras, {II}: {$C^*$}-algebras with rational generalized tracial rank one},
    journal = {C. R. Math. Rep. Acad. Sci. Canada},
    volume = {42},
    number = {4},
    pages = {451--539},
    year = {2020}
}

@Book{Anantharaman-DelarocheRenault,
  title      = {Amenable groupoids},
  publisher  = {L'Enseignement Math\'ematique},
  year       = {2000},
  author     = {Anantharaman-Delaroche, C. and Renault, J.},
  volume     = {36},
  series     = {Monographies de L'Enseignement Math\'ematique [Monographs of L'Enseignement Math\'ematique]},
  address    = {Geneva},
  isbn       = {2-940264-01-5},
  note       = {With a foreword by Georges Skandalis and Appendix B by E. Germain},
  mrclass    = {22A22 (22D25 43A07 46L05 46L10 46L80)},
  mrnumber   = {1799683 (2001m:22005)},
  mrreviewer = {Robert S. Doran},
  pages      = {196},
}

@Unpublished{STW2025,
  Title                    = {Nuclear  $C^*$-algebras: 99 problems},
  Author                   = { Schafhauser, Christopher and Tikuisis, Aaron and White, Stuart},
  Note                     = {To appear in Special Issue of \textit{Munster J. Math.} in memory of Eberhard Kirchberg, arXiv:2025.10902},
  Year                     = {2025}
}

@book{Kammeyer2019Introductiona,
    address = {Cham},
    series = {Lecture {Notes} in {Mathematics}},
    title = {Introduction to $\ell^2$-invariants},
    volume = {2247},
    isbn = {978-3-030-28296-7 978-3-030-28297-4},
    url = {http://link.springer.com/10.1007/978-3-030-28297-4},
    doi = {10.1007/978-3-030-28297-4},
    language = {en},
    urldate = {2023-02-20},
    publisher = {Springer International Publishing},
    author = {Kammeyer, Holger},
    year = {2019},
}

@article{CheegerGromov1986L2cohomology,
    title = {{$L_2$}-cohomology and group cohomology},
    volume = {25},
    url = {https://www.sciencedirect.com/science/article/pii/004093838690039X},
    number = {2},
    urldate = {2026-07-02},
    journal = {Topology},
    publisher = {Elsevier},
    author = {Cheeger, Jeff and Gromov, Mikhael},
    year = {1986},
    pages = {189--215},
}

\end{document}